\documentclass[letterpaper,11pt]{article}
\usepackage{times}
\usepackage{empheq}
\usepackage{amsmath, amssymb}
\usepackage{amsthm}
\usepackage[T1]{fontenc}
\usepackage[english]{babel}
\usepackage[utf8]{inputenc}
\usepackage{graphicx}
\usepackage{hyperref}
\usepackage{mathrsfs}
\usepackage{mathabx}
\usepackage{color}
\usepackage{anysize}
\usepackage{apacite}
\usepackage{natbib}
\usepackage{xcolor}
\marginsize{2.5cm}{2.5cm}{2.5cm}{2.5cm}
\usepackage[toc,page]{appendix}
\usepackage{url}
\usepackage{stmaryrd}
\usepackage{bbm}
\usepackage{mathtools}
\usepackage{subfigure}
\usepackage[font={small,sc}]{caption}
\usepackage{epigraph}
\usepackage{ragged2e}

\newtheorem{Theo}{Theorem}[section]
\newtheorem{Def}{Definition}[section]
\newtheorem{Prop}{Proposition}[section]
\newtheorem{Corol}{Corollary}[section]
\newtheorem{Lemma}{Lemma}[section]
\theoremstyle{definition}

\usepackage{upgreek}

\usepackage{setspace}

\AtBeginDocument{}

\DeclareMathOperator\supp{supp}

\newcommand{\Rd}{\mathbb{R}^{d}}
\newcommand{\Rm}{\mathbb{R}^{m}}
\newcommand{\RdxR}{\Rd\times\mathbb{R}}

\newcommand{\RdxRm}{\Rd\times\Rm}
\newcommand{\BorelRd}{\mathcal{B}(\mathbb{R}^{d})}

\newcommand{\BoundedBorelRd}{\mathcal{B}_{B}(\mathbb{R}^{d})}

\newcommand{\SchwartzRd}{\mathscr{S}(\Rd)}

\newcommand{\TemperedRd}{\mathscr{S}'(\Rd)}

\newcommand{\OMRd}{\mathcal{O}_{M}(\Rd)}
\newcommand{\Rp}{\mathbb{R}^{+}}
\newcommand{\RdxRp}{\mathbb{R}^{d}\times\mathbb{R}^{+}}

\newcommand{\CfdRdRp}{C_{FD}( \RdxRp )}

\newcommand{\MsgRp}{\mathscr{M}_{SG}(\Rp)}

\newcommand{\SchwartzRdxR}{\mathscr{S}(\RdxR)}
\newcommand{\TemperedRdxR}{\mathscr{S}'(\RdxR)}

\newcommand{\RdxRneg}{\Rd\times\mathbb{R}^{-}_{*}}
\newcommand{\VRdCfdRp}{\mathscr{V}(\Rd)\widehat{\boxtimes}C_{FD}(\Rp)}
\newcommand{\VRdCfdRpDual}{\left(\mathscr{V}(\Rd)\widehat{\boxtimes}C_{FD}(\Rp)\right)'}
\newcommand{\VRd}{\mathscr{V}(\Rd)}
\newcommand{\VRdDual}{\mathscr{V}'(\Rd)}
\newcommand{\VRdxRDual}{\mathscr{V}'(\RdxR)}

\newcommand{\R}{\mathbb{R}}
\newcommand{\Rneg}{\mathbb{R}^{-}_{*}}
\newcommand{\MRd}{\mathscr{M}(\Rd)}
\newcommand{\MfRd}{\mathscr{M}_{F}(\Rd)}
\newcommand{\MsgRd}{\mathscr{M}_{SG}(\Rd)}
\newcommand{\CfdRd}{C_{FD}(\Rd)}
\newcommand{\CfdRdRpDual}{C_{FD}'( \RdxRp )}
\newcommand{\CfdR}{C_{FD}(\R)}
\newcommand{\MsgRdxRp}{\mathscr{M}_{SG}(\RdxRp)}
\newcommand{\CfdRdxR}{C_{FD}(\RdxR)}
\newcommand{\CfdRdxRDual}{C_{FD}'(\RdxR)}
\newcommand{\MsgRdxR}{\mathscr{M}_{SG}(\RdxR)}
\newcommand{\OMR}{\mathcal{O}_{M}(\R)}
\newcommand{\VRdCfdR}{\mathscr{V}(\Rd)\widehat{\boxtimes}\CfdR}
\newcommand{\VRdCfdRDual}{(\mathscr{V}(\Rd)\widehat{\boxtimes}\CfdR)'}
\newcommand{\CfdRdxRp}{C_{FD}(\RdxRp)}
\newcommand{\DRd}{\mathscr{D}(\Rd)}
\newcommand{\DR}{\mathscr{D}(\R)}

\newcommand{\MR}{\mathscr{M}(\R)}
\newcommand{\DistributionsRd}{\mathscr{D}'(\Rd)}
\newcommand{\DistributionsR}{\mathscr{D}'(\R)}

\begin{document}

\begin{Large}
\begin{center}
\textbf{First-order linear evolution equations with càdlàg-in-time solutions}
\end{center}
\end{Large}
\begin{center}
Ricardo Carrizo Vergara
\end{center}
\begin{center}
June 2019
\end{center}

\begin{center}
\textbf{Abstract}
\end{center}

\small
\justify
In this work we study first-order linear parabolic evolution PDEs over $\RdxR$ and $\RdxRp$ comprising a spatial operator defined through a symbol function and a source term such that its spatial Fourier transform is a slow-growing measure over $\RdxR$. When the source term is required to has its support on $\RdxRp$, it is shown that there exists a unique solution such that its spatial Fourier transform is a slow-growing measure with support in $\RdxRp$, which in addition has a càdlàg-in-time behaviour. This allows to well-pose and analyse an initial value problem associated to this class of equations and to consider cases where the spatial operator can be a pseudo-differential operator. We also look at for solutions to the cases where the source term is  such that its spatial and spatio-temporal Fourier transforms are slow-growing measures over $\RdxR$. In such a case, it is shown that when the real part of the symbol function of the spatial operator is inferiorly bounded by a strictly positive constant, there exists a unique solution whose both spatial and spatio-temporal Fourier transforms are slow-growing measures over $\RdxR$, which also has a càdlàg-in-time behaviour. In addition, it is proven that the solution to an associated Cauchy problem converges spatio-temporally asymptotically to this unique solution as the time flows long enough.

\bigskip


\normalsize

\section{Introduction}
\label{Sec:Introduction}

In this work we look at for solutions with a càdlàg-in-time behaviour to PDEs ("\textit{P}" stands for \textit{Partial} or \textit{Pseudo}, as pleasure) over the Euclidean space-time of the form
\begin{equation}
\label{Eq:ThePDE}
\frac{\partial U}{\partial t} + \mathcal{L}_{g}U = X,
\end{equation}
where the \textit{source term} $X$ is such that its spatial Fourier transform is a slow-growing measure over $\RdxR$ and $\mathcal{L}_{g}$ is a purely-spatial pseudo-differential operator defined through the spatial Fourier transform and a symbol function $g : \Rd \to \mathbb{C}$ through $\mathcal{L}_{g} = \mathscr{F}_{S}^{-1}(g\mathscr{F}_{S}( \cdot ) )$. We consider such kinds of source terms since these act as measures in the temporal component, so it is expected that the potential solutions can be interpreted as a càdlàg function in time. Hence, an initial value problem associated to Eq. \eqref{Eq:ThePDE} make sense, even if the solution is considered, in principle, to belong to a space of distributions which have not necessarily point-wise evaluation meaning. We study such equations in a parabolic framework, for which we mean that the symbol function $g$ is such that its real part, $g_{R}$, satisfies $g_{R} \geq 0 $. Our analysis is restricted to $\RdxR$ or $\RdxRp$, without concerning about any particular spatial boundary condition.

This kind of equation is used to model a big variety of phenomena, and it can often be found in physics, biology, fluid mechanics and other fields. Linear versions of the heat or diffusion equation, some transport equations, advection-diffusion equations with damping, some forms of the Schrödinger equation and fractional (in space) versions of these equations are of the form \eqref{Eq:ThePDE}. There is a huge literature on equations of this form and some variants of these. We refer for example to classical results exposed in \citet[][Chapter XXIII]{hormander2007analysis}, \citet{lions1963quelques} , \citet{dautray2012mathematical} and the references therein. More recent works which work out non-linear and non-local forms of such equation can be found for example in \citep{avalishvili2011nonclassical}, \citep[][Chapter 30]{zeidler2013nonlinear} and in the references therein. The main general setting for the analysis is in the context of Hilbert or Banach spaces in space, while using either continuous functions in time or a distributional sense in time. For instance, the classical Lions-Tartar theorem and some of its generalizations \citep{simon2007generalizacion} are under such a setting.

In this work we study this equation in a particularly different setting from what can be found in the literature. We do not focus on the theory of Hilbert spaces. We analyse equation \eqref{Eq:ThePDE} in a distributional sense, for which we select a suitable space of spatio-temporal distributions where $X$ belongs. Namely, we require $X$ to belong to the space of tempered distributions over $\RdxR$ whose spatial Fourier transforms are slow-growing measures over $\RdxR$, here denoted $\VRdCfdRDual$. Such space of distributions involves a lot of classically used Hilbert and Banach spaces, for instance the spaces $L^{p}(\RdxR)$ for $p \in \left[1,2\right]$ (and hence some Sobolev spaces) and the space of finite complex (Radon) measures over $\RdxR$. We show that this requirement on $X$ implies that the solution, when it exists, has a càdlàg-in-time behaviour since it first temporal derivative acts as a measure in time. In addition, it is proven that the solution $U$ also belongs to this space of distributions. The fact that the solution $U$ is also a member of $\VRdCfdRDual$ allows us to gain freedom on the selection of the symbol function $g$ with respect to the case of general tempered distributions, where $g$ is needed to be in the space $\OMRd$ (see Sections \ref{Sec:Preliminaries} and \ref{Sec:OperatorsSymbol} for more details). For instance, when $g$ is a continuous and polynomially bounded function the operator $\mathcal{L}_{g}$ can be freely applied over any distribution in $\VRdCfdRDual$. Hence, some fractional pseudo-differential operators such as $(-\Delta)^{\frac{\alpha}{2}}$ for $\alpha \geq 0 $ are involved in our scope. The analysis we present in this work has been already presented in broad terms in \citet[][Appendix C]{carrizov2019development}. Here we present it with more depth and in a shorter and self-contained exposition.

Another motivation to consider cases where $X$ belongs to a suitable space of distributions which act as measures in time but not necessarily as members of a Hilbert space comes, actually, from the world of Stochastic Analysis. The main objective of this paper is tu set-up the deterministic analogous framework for the analysis of a \textit{stochastized} version of the PDE \eqref{Eq:ThePDE} (hence, a SPDE) in which the source term $X$ may represent a spatio-temporal White Noise, a coloured-in-space and white-in-time noise, or, more generally, a stationary generalized random field with a measure-in-time behaviour. Such mathematical objects cannot be always interpreted as members of a Hilbert space, but they do can be interpreted as generalized stochastic processes whose spatial Fourier transforms are random measures in a mean-square sense. Such analysis has been worked out in \citet[][Appendix C]{carrizov2019development} and it will be exposed in depth in a forthcoming paper.

This paper is organised as follows. In Section \ref{Sec:TheoreticalBackground} we set-up the tools we will use and the main mathematical results for those. We make reminders on measures, distributions and topological vector spaces in Section \ref{Sec:Preliminaries}. In Section \ref{Sec:SlowGrowingMeasures} we introduce the space of slow-growing complex measures $\MsgRd$, which can be identified as the dual of the space of continuous fast-decreasing functions, $\CfdRd$. We introduce different spaces of test-functions used to work in spatio-temporal frameworks in Section \ref{Sec:ConvinientTestFunctions}. We introduce the space $\VRd$ of functions which are Fourier transforms of functions in $\CfdRd$, together with its dual space $\VRdDual$ of tempered distributions whose Fourier transforms are in $\MsgRd$. We introduce the associate space $\VRdCfdR$ of test-functions which are spatial Fourier transforms of functions in $\CfdRdxR$, and its dual space $\VRdCfdRDual$ of spatio-temporal distributions whose spatial Fourier transforms are in $\MsgRdxR$. We introduce the spaces of restrictions to $\RdxRp$ of members of the aforementioned spaces. In Section \ref{Sec:OperatorsSymbol} we define properly operators of the form $\mathcal{L}_{g}$ and we show how to solve PDEs associated to those under suitable conditions. In Section \ref{Sec:CadlagInTimeRepresentations} we introduce our concept of distribution having a càdlàg-in-time behaviour, and we give a Theorem which presents sufficient conditions for a distribution belonging to $\VRdCfdRDual$ or to $\MsgRdxR$ to have a càdlàg-in-time representation. Section \ref{Sec:AnalysisPDE} is devoted to the analysis of the PDE \eqref{Eq:ThePDE} without yet explicitly concerning about a Cauchy problem associated to it. In Section \ref{Sec:SolutionsOverRdxRp} we look at for solutions over $\RdxRp$ requiring the source term to has its support on $\RdxRp$. We show that there exists a unique solution in $\VRdCfdRDual$ with support in $\RdxRp$ and that it has a càdlàg-in-time representation which we make explicit. In Section \ref{Sec:SolutionsOverRdxR} we look at for solutions over $\RdxR$ when the source term is in $\VRdxRDual \cap \VRdCfdRDual$. When requiring the symbol function to satisfy $g_{R} \geq \kappa > 0 $, we prove that there exists a unique solution in  $\VRdCfdRDual\cap\VRdxRDual$ which has a càglàg-in-time representation which we make explicit. In Section \ref{Sec:CauchyProblem} we analyse a Cauchy problem associated to \eqref{Eq:ThePDE} over $\RdxRp$ with initial conditions in $\VRdDual$ for the case $g_{R} \geq 0$. We obtain suitable existence and uniqueness results together with the càdlàg-in-time representation of the solution. In Section \ref{Sec:LongTermAsymptotics} we show the relationship between the solution of a Cauchy problem obtained in Section \ref{Sec:CauchyProblem} and the solution in $\VRdxRDual$ obtained in Section \ref{Sec:SolutionsOverRdxR} for the case $g_{R} \geq \kappa
 > 0 $. We prove that the solution to the Cauchy problem converges \textit{spatio-temporally} to the solution in $\VRdxRDual$ as the time flows long enough. We also show that if the initial condition in the Cauchy problem is the evaluation at $0$ of the unique solution in $\VRdxRDual$, then the solution to the Cauchy problem and the unique solution in $\VRdxRDual$ coincide over $\RdxRp$.

We set up some \textit{almost-french} notations. We denote $\mathbb{N} := \lbrace 0 , 1 , 2 , ... \rbrace$, $\mathbb{N}_{*} = \mathbb{N}\setminus \lbrace 0 \rbrace $, $\Rp = \left[ 0 , \infty \right)$, $\R^{-} = \left( -\infty , 0 \right]$, $\Rp_{*} = \Rp\setminus\lbrace 0 \rbrace$, $\Rneg = \R^{-}\setminus\lbrace 0 \rbrace$. $\mathbf{1}_{A}$ denotes the indicator function of the set $A$. For $x \in \Rd$, $|x|$ denotes the Euclidean norm of $x$, regardless of the dimension $d$.

\section{Theoretical Background}
\label{Sec:TheoreticalBackground}

\subsection{Preliminaries}
\label{Sec:Preliminaries}

Let us make precise the definition of complex Radon measure over $\Rd$ we use, with $d \in \mathbb{N}_{*}$. We denote $\BorelRd$ the collection of Borel subsets of $\Rd$ and $\BoundedBorelRd$ the collection of bounded Borel subsets of $\Rd$. 

\begin{Def}
\label{Def:ComplexRadonMeasure}
A complex Radon measure over $\Rd$ is a function $\mu : \BoundedBorelRd \to \mathbb{C}$ such that for every countable collection of mutually disjoint bounded Borel sets $(A_{n})_{n \in \mathbb{N}} \subset \BoundedBorelRd$ such that $\bigcup_{n \in \mathbb{N}}A_{n} \in \BoundedBorelRd$, the $\sigma-$additivity property holds:
\begin{equation}
\label{Eq:SigmaAdditivity}
\mu\left( \bigcup_{n \in \mathbb{N} }A_{n} \right) = \sum_{n \in \mathbb{N}} \mu(A_{n}).
\end{equation}
\end{Def}

We note $\MRd$ the set of complex Radon measures over $\Rd$. $\MRd$ is a complex vector space. For a measure $\mu \in \MRd$, its total variation measure $|\mu|$ is defined as
\begin{equation}
\label{Eq:DefTotalvariationMeasure}
|\mu|(A) := \sup \bigg\lbrace  \sum_{n \in \mathbb{N}}  |\mu(E_{n})| \ \big| \  (E_{n})_{n \in \mathbb{N}} \subset \BorelRd  \ \mathrm{partition \ of } \ A   \bigg\rbrace, \quad \forall A \in \BoundedBorelRd.
\end{equation}
$|\mu|$ is, indeed, a (Radon) measure \citep[Theorem 6.2, applicable over bounded subsets]{rudin1987real}, which is positive. A measure $\mu \in \MRd$ is said to be finite if $\sup_{A \in \BoundedBorelRd} |\mu|(A) < \infty. $ In such a case, the value $\mu(A)$ can be well-defined for every $A \in \BorelRd$ \citep[Theorem 6.4]{rudin1987real}. We note $\MfRd$ the space of finite measures in $\MRd$. If $\mu \in \MRd$, its support is defined as the complementary of the biggest open set where the total variation of $\mu$ is null:
\begin{equation}
\label{Eq:DefSupportMeasure}
\supp(\mu)  := \left( \bigcup \lbrace O \subset \mathbb{R}^{d} \ : \ O \mathrm{\ is \ open \ and  } \ |\mu|(O) = 0   \rbrace \right)^{c}. 
\end{equation} 
For a measurable and locally bounded function $f : \Rd \mapsto \mathbb{C}$, and for $\mu \in \MRd$, we denote $f\mu$ the multiplication measure $f\mu : A \in \BoundedBorelRd \mapsto  \int_{A}f(x)d\mu(x)$, which is indeed in $\MRd$. In particular, for any $A \in \BorelRd$, the measure $\mathbf{1}_{A}\mu$ is the \textit{restriction measure} of $\mu$ over $A$.

We denote $\DRd$ the space of compactly supported complex smooth functions over $\Rd$ and $\DistributionsRd$ its dual space of distributions over $\Rd$.

For the case $d=1$, the following result is widely known: for any $\mu \in \MR$, there exists a unique càdlàg function (that is, right-continuous with left-limits) $f: \R \to \mathbb{C}$ such that $\frac{df}{dt} = \mu $ in distributional sense and such that $f(0)=0$. Such a function is explicitly given by 
\begin{equation}
\label{Eq:CadlagPrimitiveMeasure1D}
f(t) = \mu(\left( 0 , t \right])\mathbf{1}_{\Rp}(t) - \mu(\left( t , 0 \right)  )\mathbf{1}_{\Rneg}(t), \quad \forall t \in \R.
\end{equation}
Conversely, if a distribution $T \in \DistributionsR$ is such that $\frac{dT}{dt} \in \MR$, then $T$ can be represented by a unique càdlàg function, for which it holds
\begin{equation}
\label{Eq:DifferenceCadlagMeasure}
T(t) - T(s) = \frac{dT}{dt}(\left( s , t \right]  ), \quad \forall t,s \in \R, t>s.
\end{equation}

For any $d \in \mathbb{N}_{*}$, we denote $\SchwartzRd$ the Schwartz space of smooth functions with fast-decreasing behaviour, and $\TemperedRd$ its dual space of tempered distributions over $\Rd$. We use the following convention for the Fourier transform $\mathscr{F} : \SchwartzRd \to \SchwartzRd$ and its inverse:
\begin{equation}
\label{Eq:FourierTransformAndInverse}
\mathscr{F}(\varphi)(\xi) = \frac{1}{(2\pi)^{\frac{d}{2}}}\int_{\Rd}e^{-i\xi^{T}x}\varphi(x) dx \quad ; \quad \mathscr{F}^{-1}(\varphi)(\xi) = \frac{1}{(2\pi)^{\frac{d}{2}}}\int_{\Rd}e^{i\xi^{T}x}\varphi(x) dx.
\end{equation}
The Fourier transform over $\TemperedRd$ is the adjoint of the Fourier transform over $\SchwartzRd$. Analogously for its inverse. We denote $\OMRd$ the space of multiplicators of the Schwartz space (complex smooth functions with polynomially bounded derivatives of all orders), and $\mathcal{O}_{c}'(\Rd)$ the space of convoluters of tempered distributions (\textit{fast-decrasing distributions}, see \citet[Chapter VII, \S 5]{schwartz1966theorie} or \citet[Chapter 30]{treves1967topological}). It holds $\mathscr{F}(\OMRd) = \mathcal{O}_{c}'(\Rd)$. The exchange formula of the Fourier transform with our convention is given by
\begin{equation}
\label{Eq:ExchangeFormula}
\mathscr{F}( S \ast T) = (2\pi)^{\frac{d}{2}}\mathscr{F}(S) \mathscr{F}(T)\quad ; \quad \mathscr{F}(Tg) =  (2\pi)^{-\frac{d}{2}} \mathscr{F}(T) \ast \mathscr{F}(g),
\end{equation}
for every  $T \in \TemperedRd$, $g \in \OMRd$, and $ S \in \mathcal{O}_{c}'(\Rd)$.
   
We make some reminders on topological vector spaces and linear operators between them. A (complex) Hausdorff locally convex topological vector space (HLCTVS) $E$ is a complex vector space equipped with a family of semi-norms $(p_{j})_{j \in J}$ with the axiom $(p_{j}(x) = 0 \hbox{ for all } j \in J  ) \Rightarrow x = 0  $. The topology induced on $E$ is the weakest topology which makes continuous the semi-norms $(p_{j})_{j \in J}$, the addition and the scalar multiplication. $(p_{j})_{j \in J}$ is said to be directed if for every $j_{1},j_{2} \in J$ there exists $j_{3} \in J$ and $C > 0 $ such that $p_{j_{1}}(x) + p_{j_{2}}(x) \leq C p_{j_{3}}(x)$ for all $x \in E$. We can always find a family of directed semi-norms which determines the topology of a HLCTVS. If the topology of $E$ is determined by a countable directed family of semi-norms, then its topology is metric, and if it is in addition complete, $E$ is called a Fréchet space. Let us suppose de family $(p_{j})_{j \in J}$ is directed. If $F$ is another HLCTVS with topology determined by a directed family of semi-norms $(q_{k})_{k \in K}$, then a linear function $\mathcal{L} : E \to F $ is continuous if and only if for every $k \in K$ there exist $C \geq 0 $ and $j \in J$ such that  $q_{k}(x) \leq C p_{j}(x)$ for all $x \in E$. If $E$ is a complex HLCTVS, we denote $E'$ its dual space, that is, the space of continuous linear functions from  $E$ to $\mathbb{C}$.  If $T \in E'$ and $x \in E$, we denote $\langle T , x \rangle := T(x) \in \mathbb{C}$. We endow $E'$ with the weak-$*$ topology, which is the topology induced by the family of semi-norms $(q_{x})_{x \in E}$ given by $q_{x}(T) = |\langle T , x \rangle|$, for all $x \in E$. $E'$ with this topology is a complex HLCTVS.  If $E$ and $F$ are two HLCTVSs, and if $\mathcal{L} : E \to F $ is a continuous linear operator, then its \textit{adjoint} is the operator $\mathcal{L}^{*} : F' \to E'$ determined by $ \langle \mathcal{L}^{*}(T) , x \rangle = \langle T , \mathcal{L}(x) \rangle$, for every $ T \in F'$ and every $x \in E$. $\mathcal{L}^{*}$ is a continuous linear operator from $F'$ to $E'$. See \citet[Chapter V]{reed1980methods}  for more details.

\subsection{Slow-growing measures}
\label{Sec:SlowGrowingMeasures}

\begin{Def}
\label{Def:SlowGrowingMeasure}
Let $\mu \in \MRd$. We say that $\mu$ is slow-growing if there exists a strictly positive polynomial $p : \Rd \to \Rp_{*}$ such that the measure $\frac{d\mu(x)}{p(x)}$ is finite, or equivalently, if there exists $N \in \mathbb{N}$ such that $\int_{\Rd}\frac{d|\mu|(x)}{(1+|x|^{2})^{N}} < \infty.$ 
\end{Def}

It is immediate that a measure $\mu \in \MRd$ is slow-growing if and only if its total variation measure also is it. We denote $\MsgRd$ the complex vector space of slow-growing measures over $\Rd$. It is clear that if $f : \Rd \to \mathbb{C}$ is a measurable and polynomially bounded function and $\mu \in \MsgRd$, then the multiplication measure $f\mu$ is in $\MsgRd$. Every slow-growing measure $\mu$ is identified with a unique tempered distribution through $ \varphi \in \SchwartzRd \mapsto \langle \mu , \varphi \rangle := \int_{\Rd}\varphi(x)d\mu(x) $, having hence $\MsgRd \subset \TemperedRd$.

We denote $C(\Rd)$ the space of continuous complex functions over $\Rd$. Let us introduce the space of fast-decreasing continuous functions:
\begin{equation}
\label{Def:CfdRd}
\CfdRd := \lbrace \varphi \in C(\Rd) \ \big| \ \sup_{x \in \Rd}\left| (1+|x|^{2})^{N}\varphi(x)  \right| < \infty, \ \forall N \in \mathbb{N} \rbrace. 
\end{equation}
Obviously $\SchwartzRd \subset \CfdRd$. We endow $\CfdRd$ with the topology induced by the directed family of semi-norms given by $p_{N}(\varphi) = \sup_{x \in \Rd}| (1+|x|^{2})^{N}\varphi(x)|$, for every $N \in \mathbb{N}$. $\CfdRd$ is a Fréchet space; its completeness can be proven using the same classical arguments which prove the completeness of $\SchwartzRd$. A linear functional $T : \CfdRd \to \mathbb{C}$ is continuous if and only if there exists $C > 0 $ and $N \in \mathbb{N}$ such that
\begin{equation}
\label{Eq:ConditionContinuousLinearFunctionalCfd}
\left| \langle T , \varphi \rangle \right| \leq C \sup_{x \in \Rd} \left| (1+|x|^{2})^{N}\varphi(x) \right|, \quad \forall \varphi \in \CfdRd.
\end{equation}

We denote by $C_{c}(\Rd)$ the space of compactly supported continuous complex functions over $\Rd$, and by $C_{0}(\Rd)$ the space of continuous complex functions over $\Rd$ vanishing at infinity. We recall that for these spaces, when endowed with suitable topologies, the Riesz Representation Theorem holds, for which we have $C_{c}'(\Rd) = \MRd$ \citep[Section V.4, Example 1]{reed1980methods} and $C_{0}'(\Rd) = \MfRd$ \citep[Theorem 6.19]{rudin1987real}. We present here the case of $\CfdRd$ and $\MsgRd$.

\begin{Theo}[\textbf{Riesz Representation for slow-growing measures}]
\label{Theo:RieszReperesentationDualCfdMsg}
$\mathscr{M}_{SG}(\Rd) = C_{FD}'(\Rd)$, that is, every slow-growing measure $\mu \in \mathscr{M}_{SG}(\Rd)$ defines a continuous linear functional $T$ over $C_{FD}(\Rd)$ through the integral
\begin{equation}
\label{Eq:RepresentationLinearFunctionalCfdMsg}
\langle T , \varphi \rangle = \int_{\Rd}\varphi(x)d\mu(x), \quad \forall \varphi \in C_{FD}(\Rd).
\end{equation}
Conversely, for every continuous linear functional $T:C_{FD}(\Rd) \to \mathbb{C}$ there exists a unique $\mu \in \mathscr{M}_{SG}(\Rd)$ such that (\ref{Eq:RepresentationLinearFunctionalCfdMsg}) holds.
\end{Theo}

\textbf{Proof: } The proof of $\MsgRd \subset C_{FD}'(\Rd)$ is left to the reader. Let us prove that $C_{FD}'(\Rd) \subset \MsgRd$. Let $T \in C_{FD}'(\mathbb{R}^{d})$ and let  $C >0 $ and $N \in \mathbb{N}$ such that \eqref{Eq:ConditionContinuousLinearFunctionalCfd} holds. Let us define the linear functional $(1+|x|^{2})^{-N}T : C_{FD}(\mathbb{R}^{d}) \to \mathbb{C}$ by $(1+|x|^{2})^{-N}T(\varphi) := \langle T , (1+|x|^{2})^{-N}\varphi \rangle$. Since for all $ \varphi \in C_{FD}(\mathbb{R}^{d})$ the function $ x \in \Rd \mapsto (1+|x|^{2})^{-N}\varphi(x)$ is also in $C_{FD}(\mathbb{R}^{d})$, this functional is well-defined. In addition,
\begin{equation}
\label{Eq:(1+|x|^(2))^(-N)TVarphiBounded}
|(1+|x|^{2})^{-N}T(\varphi)| = | \langle T , (1+|x|^{2})^{-N} \varphi \rangle | \leq C \sup_{x \in \Rd}\left| (1+|x|^{2})^{N}(1+|x|^{2})^{-N} \varphi(x) \right| = C \sup_{x \in \Rd}\left| \varphi(x) \right|,
\end{equation}
hence, $T$ is continuous. Expression \eqref{Eq:(1+|x|^(2))^(-N)TVarphiBounded} holds in particular for every $ \varphi \in C_{c}(\Rd) \subset C_{FD}(\Rd)$. The functional $(1+|x|^{2})^{-N}T$ is hence a bounded linear functional over $C_{c}(\Rd)$ endowed with the supremum norm. By Hahn-Banach extension Theorem \citep[Theorem III.5 or Theorem V.3]{reed1980methods}, $(1+|x|^{2})^{-N}T$ can be extended to a continuous linear functional over $C_{0}(\Rd)$, and since $C_{c}(\Rd)$ is dense in $C_{0}(\Rd)$ the extension is unique. Eq. \eqref{Eq:(1+|x|^(2))^(-N)TVarphiBounded} holds for every $\varphi \in C_{0}(\Rd)$. Since $C_{0}'(\Rd) = \mathscr{M}_{F}(\Rd)$, we conclude that $(1+|x|^{2})^{-N}T$ is identified with a unique measure $\nu \in \mathscr{M}_{F}(\mathbb{\Rd})$. Consider then the multiplication measure $d\mu(x) = (1+|x|^{2})^{N}d\nu(x) $, which is in $\mathscr{M}_{SG}(\Rd)$. We conclude that for every $\varphi \in C_{FD}(\Rd)$ we have
\small
\begin{equation}
\label{Eq:ConclusionRRTMsg}
\langle T , \varphi \rangle = \langle T , \frac{(1+|x|^{2})^{N}}{(1+|x|^{2})^{N}}\varphi \rangle = \langle (1+|x|^{2})^{-N}T , (1+|x|^{2})^{N}\varphi \rangle = \int_{\Rd}\varphi(x)(1+|x|^{2})^{N}d\nu(x) = \int_{\Rd}\varphi(x)d\mu(x). \blacksquare
\end{equation}
\normalsize

\subsection{Convenient spaces of test-functions and their dual spaces}
\label{Sec:ConvinientTestFunctions}

Let us denote $C_{0}^{\infty}(\Rd)$ the space of smooth complex functions over $\Rd$ such that all of their derivatives vanish at infinity.  Let us introduce the next space of test-functions:
\begin{equation}
\label{Eq:DefVRd}
\VRd = \lbrace  \varphi \in C_{0}^{\infty}(\Rd) \ \big| \ \exists \phi \in C_{FD}(\Rd) \hbox{ such that } \varphi = \mathscr{F}(\phi)   \rbrace = \mathscr{F}(  C_{FD}(\Rd) ).
\end{equation}
$\mathscr{V}(\Rd)$ is defined in such a way that $\mathscr{V}(\Rd) \subset C_{0}^{\infty}(\Rd)$, but this requirement actually follows from Riemann-Lebesgue Lemma. $\mathscr{V}(\Rd)$ is a strict subspace of $C_{0}^{\infty}(\Rd)$. The space $\mathscr{V}(\Rd)$ can be defined equivalently as the space of functions in $C_{0}^{\infty}(\Rd)$ such that their inverse Fourier transforms in distributional sense can be identified with a function in $C_{FD}(\Rd)$. Of course, if we use the inverse Fourier transform instead of the Fourier transform in the definition \eqref{Eq:DefVRd} of $\mathscr{V}(\Rd)$, the space keeps unchanged.

We endow $\mathscr{V}(\Rd)$ with the topology induced by the directed family of semi-norms 
\begin{equation}
\label{Eq:SeminormsVRd}
p_{N}(\varphi) = \sup_{\xi \in \Rd}\left|(1+|\xi|^{2})^{N}\mathscr{F}(\varphi)(\xi)  \right|, \quad N \in \mathbb{N}.
\end{equation}
It is quite immediate that the Fourier transform and its inverse interpreted in distributional sense define continuous and bijective linear transformations from $\CfdRd$ to $\VRd$ (or from $\VRd$ to $\CfdRd$, as pleasure). The completeness of $\CfdRd$ and the continuity of the Fourier transform imply the completeness of $\VRd$. $\VRd$ is a Fréchet space.

\begin{Prop}
\label{Prop:SchwartzRdDenseCfdVRd}
The Schwartz space satisfies $\SchwartzRd \subset \mathscr{V}(\Rd)$ and $\SchwartzRd \subset C_{FD}(\Rd)$ and it is a dense subspace of both spaces with their respective topologies.
\end{Prop}
The next Lemma will be widely used in this work.
\begin{Lemma}
\label{Lemma:ConvexityInequalities}
For every $x,y \in \Rd$ and for every $m \in \mathbb{N}$, $(1+|x|^{2})^{m} \leq 2^{m-1}\left[ (1+2|x-y|^{2})^{m} + 2^{m}|y|^{2m}  \right]$. In particular, $(1+|x|^{2})^{m} \leq 2^{m-1}\left( 3^{m} + 2^{m}|y|^{2m}  \right)$ when $|x-y|\leq 1$.
\end{Lemma}
\textbf{Proof: } Using first the convexity of the function $ x \in \Rd \mapsto 1+|x|^{2}$ and then the convexity of the function $ t \in \Rp \mapsto t^{m} $, we have
\begin{equation}
\label{Eq:ProofConvexityInequalities}
\left. \begin{aligned}
(1+|x|^{2})^{m} &= \left( 1+|x-y + y|^{2}  \right)^{m} = \left( 1+\left| \frac{2(x - y)}{2} + \frac{2y}{2} \right|^{2}  \right)^{m} \leq \left(1 + 2|x-y|^{2} + |y|^{2} \right)^{m} \\
&= \left( \frac{2(1+2|x-y|^{2})}{2} + \frac{4|y|^{2}}{2}  \right)^{m} \leq \frac{\left[ 2(1+2|x-y|^{2}) \right]^{m}}{2} + \frac{(4|y|^{2})^{m}}{2} \\ 
&= 2^{m-1}\left[ (1+2|x-y|^{2})^{m} + 2^{m}|y|^{2m}  \right]. \quad \blacksquare
\end{aligned} \right.
\end{equation}

\textbf{Proof of Proposition  \ref{Prop:SchwartzRdDenseCfdVRd}:} The inclusions are straightforward. We will just prove the density of $\SchwartzRd$ in $C_{FD}(\Rd)$. The density in $\mathscr{V}(\Rd)$ follows immediately from the continuity of the Fourier transform. 

We first prove that if $f \in C_{FD}(\mathbb{R}^{d})$ and $\varphi \in \SchwartzRd$, then $ f \ast \varphi \in \SchwartzRd$. It is clear that $f$ is integrable and bounded, as well as $\varphi$ which is in addition smooth. Thus $f \ast \varphi $ is a smooth integrable and bounded function, and its Fourier transform satisfies $\mathscr{F}( f \ast \varphi) = (2\pi)^{\frac{d}{2}}\mathscr{F}(f)\mathscr{F}(\varphi)$. Since $f \in C_{FD}(\mathbb{R}^{d})$, then $\mathscr{F}(f) \in \mathscr{V}(\Rd) \subset C_{0}^{\infty}(\Rd) \subset \OMRd$. This implies that $(2\pi)^{\frac{d}{2}}\mathscr{F}(f)\mathscr{F}(\varphi) \in \SchwartzRd$. This proves that $f \ast \varphi = \mathscr{F}^{-1}\left(  (2\pi)^{\frac{d}{2}}\mathscr{F}(f)\mathscr{F}(\varphi) \right) \in \SchwartzRd$.

Let $(\phi_{n})_{n \in \mathbb{N}} \subset  \DRd$ be a regularizing sequence of positive functions, such that $\int_{\mathbb{R}^{d}}\phi_{n}(x)dx = 1 $  and $\supp(\phi_{n}) = \overline{B_{\frac{1}{n}}(0)}$ for all $n\in\mathbb{N}$, where $B_{r}(0)$ denotes the open ball of radius $r > 0 $ centred at $0$. We consider the sequence of functions $f_{n} = f \ast \phi_{n}$, which are all in $\SchwartzRd$. We will prove that $f_{n} \to f $ in $\CfdRd$. Let $m \in \mathbb{N}$ be fixed.  Let $\epsilon > 0 $. As $f \in  C_{FD}(\mathbb{R}^{d})$, we can take $R > 0 $ large enough such that for every $x$ such that $|x| > R-1$, $(1+2|x|^{2})^{m}|f(x)| < \frac{\epsilon}{3(2^{m-1}+2^{2m-1})} $. Notice that in this case, $(1+|x|^{2})^{m}|f(x)| < \frac{\epsilon}{3} $. Since $f$ is continuous, it is uniformly continuous over the compact set $\overline{B_{R+1}(0)}$.  Thus, there exists $\delta > 0 $ such that if $|x-y|<\delta$, then $|f(x)-f(y)|< \frac{\epsilon}{3(1+R^{2})^{m}}$ for all $x,y \in \overline{B_{R+1}(0)}$. Consider $n_{0} \in \mathbb{N}$ such that $\frac{1}{n_{0}} < \delta $. Then,  for all $n \geq n_{0}$,
\small
\begin{equation}
\label{Eq:BoundingUnifConvDensitySRdCfdRd}
\left.\begin{aligned}
\sup_{x \in \Rd}\left| (1+|x|^{2})^{m}(f(x) - f_{n}(x)) \right|  &=  \sup_{x \in \mathbb{R}^{d}}\left| \int_{B_{\frac{1}{n}}(0)}(1+|x|^{2})^{m}(f(x)-f(x-y))\phi_{n}(y)dy \right|   \\ 
  &\leq  \sup_{x \in \overline{B_{R}(0)}}\left| \int_{B_{\frac{1}{n}}(0)}(1+|x|^{2})^{m}(f(x)-f(x-y))\phi_{n}(y)dy \right| \\
  &\quad + \underbrace{\sup_{x \in \overline{B_{R}(0)}^{c}}\left| \int_{B_{\frac{1}{n}}(0)}(1+|x|^{2})^{m}(f(x)-f(x-y))\phi_{n}(y)dy \right|}_{(a)}. \\ 
\end{aligned}\right.
\end{equation}
\normalsize
For the first term, the uniform continuity of $f$ implies
\small
\begin{equation}
\label{Eq:BoundingPartA}
\sup_{x \in \overline{B_{R}(0)}}\left| \int_{B_{\frac{1}{n}}(0)}(1+|x|^{2})^{m}(f(x)-f(x-y))\phi_{n}(y)dy \right| \leq \int_{B_{\frac{1}{n}}(0)}(1+R^{2})^{m}\dfrac{\epsilon}{3(1+R^{2})^{m}}\phi_{n}(y)dy = \dfrac{\epsilon}{3}. 
\end{equation}
\normalsize
Regarding the second term $(a)$, the integral is split  to obtain
\begin{equation}
\label{Eq:BoundingPartB}
(a) \leq \sup_{x \in \overline{B_{R}(0)}^{c}}\Big\lbrace \underbrace{\int_{B_{\frac{1}{n}}(0)}(1+|x|^{2})^{m}|f(x)|\phi_{n}(y)dy}_{\leq \frac{\epsilon}{3} } + \int_{B_{\frac{1}{n}}(0)}(1+|x|^{2})^{m}|f(x-y)|\phi_{n}(y)dy  \Big\rbrace.
\end{equation}
Applying Lemma \ref{Lemma:ConvexityInequalities}, one has
\begin{equation}
\label{Eq:BoundingIntegralsPartB}
\left. \begin{aligned}
\int_{B_{\frac{1}{n}}(0)}(1+|x|^{2})^{m}|f(x-y)|\phi_{n}(y)dy  &\leq   2^{m-1}\Bigg[ \int_{B_{\frac{1}{n}}(0)}\underbrace{(1+2|x-y|^{2})^{m}|f(x-y)|}_{< \frac{\epsilon}{3(2^{m-1}+2^{2m-1})} \text{ from } |x-y| > R-1 }\phi_{n}(y)dy  \\
 &\quad  +  2^{m}\int_{B_{\frac{1}{n}}(0)}\underbrace{|f(x-y)|}_{ < \frac{\epsilon}{3(2^{m-1}+2^{2m-1})}}\underbrace{|y|^{2m}}_{\leq 1}\phi_{n}(y)dy \Bigg]  \\
   &<   2^{m-1}\left(\dfrac{\epsilon}{3(2^{m-1}+2^{2m-1})} + 2^{m}\dfrac{\epsilon}{3(2^{m-1}+2^{2m-1})}\right) = \dfrac{\epsilon}{3}.
\end{aligned} \right.
\end{equation}
Hence considering \eqref{Eq:BoundingPartB} and \eqref{Eq:BoundingIntegralsPartB} we finally obtain $(a) < \frac{2\epsilon}{3}$. Putting together this result and \eqref{Eq:BoundingPartA} on equation \eqref{Eq:BoundingUnifConvDensitySRdCfdRd}, we finally obtain  that for all $n \geq n_{0}$,
\begin{equation}
\sup_{x \in \Rd}\left| (1+|x|^{2})^{m}(f(x) - f_{n}(x)) \right| < \epsilon.
\end{equation}
Hence, $\sup_{x \in \Rd}\left| (1+|x|^{2})^{m}(f(x) - f_{n}(x)) \right| \to 0 $ as $n \to \infty$. Since $m$ was arbitrary, this procedure applies for all $m \in \mathbb{N}$. We therefore conclude that $f_{n} \to f $ in $\CfdRd$, hence $\SchwartzRd$ is dense in $C_{FD}(\mathbb{R}^{d})$. $\blacksquare$

Let us now consider the dual space $\VRdDual$. A linear functional $T : \mathscr{V}(\Rd) \to \mathbb{C}$ is in $\mathscr{V}'(\Rd)$ if and only if there exist $C>0$ and $N \in \mathbb{N}$ such that
\begin{equation}
\label{Eq:CriterionTinVRdDual}
|\langle T , \varphi \rangle | \leq C \sup_{\xi \in \Rd}\left| (1+|\xi|^{2})^{N}\mathscr{F}(\varphi)(\xi) \right|, \quad \forall \varphi \in \mathscr{V}(\Rd).
\end{equation}
The density of $\SchwartzRd $ in $\mathscr{V}(\Rd)$ and the integrability of the functions in $C_{FD}(\Rd)$ allow to conclude the following inclusions:
\begin{equation}
\label{Eq:InclusionsSchwartzVRd}
\SchwartzRd \subset \mathscr{V}(\Rd) \subset \mathscr{V}'(\Rd) \subset \TemperedRd.
\end{equation}
The Fourier transform over $\mathscr{V}'(\Rd)$ can be defined equivalently as the restriction over $\mathscr{V}'(\Rd)$ of the Fourier transform on $\TemperedRd$, or as the adjoint of the Fourier transform $\mathscr{F} : C_{FD}(\Rd) \to \mathscr{V}(\Rd)$. Using this second option, it follows immediately that it is a continuous  linear functional $\mathscr{F} : \mathscr{V}'(\Rd) \mapsto C_{FD}'(\Rd)$. Since $C_{FD}'(\Rd) = \mathscr{M}_{SG}(\Rd)$ from Riesz Representation Theorem \ref{Theo:RieszReperesentationDualCfdMsg}, we conclude that $\mathscr{V}'(\Rd)$ is the space of tempered distributions such that their Fourier transforms (or inverse) are slow-growing measures, having $\VRdDual = \mathscr{F}(\MsgRd)$.

Let us now introduce some special conventions and spaces of test-functions adapted to a \textit{spatio-temporal framework}. These conventions will be used all along this work unless stated otherwise. $d \in \mathbb{N}_{*}$ will always denote the spatial dimension. We will explicitly write $\RdxR$ for the space-time Euclidean space. We sometimes call $\Rp$ the \textit{positive-time}. The letters $(x,t)$ will always denote a point in the \textit{physical space-time}, and the letters $(\xi,\omega)$ will be used to denote a point in the \textit{frequency space-time}.  We will use the letters $\varphi, \phi$ for spatial test-functions, the letters $\theta, \vartheta$ for temporal test-functions and the letter $\psi$ for spatio-temporal test functions.  $\mathscr{F}$ will denote a \textit{spatio-temporal Fourier transform}, $\mathscr{F}_{S}$ a \textit{spatial Fourier transform}, and $\mathscr{F}_{T}$ a \textit{temporal Fourier transform}, all of these operations being defined over $\TemperedRdxR$. We explicit them through
\begin{equation}
\label{Eq:SpatialTemporalFourierTransforms}
\mathscr{F}_{S}(\psi)(\xi,t) = \frac{1}{(2\pi)^{\frac{d}{2}}}\int_{\Rd}e^{-i\xi^{T}x}\psi(x,t) dx \quad ; \quad \mathscr{F}_{T}(\psi)(x,\omega) = \frac{1}{\sqrt{2\pi}}\int_{\R}e^{-i\omega t}\psi(x,t) dt, 
\end{equation}
for all $\psi \in \SchwartzRdxR$. The symbol $\boxtimes$ will denote a \textit{spatio-temporal tensor product}, that is, $S \boxtimes T$ denotes the tensor product between a spatial distribution $S$ and a temporal distribution $T$. The same applies for the tensor product between spatial and temporal functions or measures. We will always use the \textit{space-time} order notation, with the spatial object \textit{at the left} and the temporal \textit{at the right} of the symbol $\boxtimes$.

We will first consider a space of test-functions acting over the whole space-time $\RdxR$:
\begin{equation}
\label{Eq:DefVRdCfdR}
\left. \begin{aligned}
\VRdCfdR  &:= \lbrace \psi \in C(\RdxR) \ \big| \ \exists \psi_{2} \in C_{FD}(\RdxR) \hbox{ such that } \psi = \mathscr{F}_{S}(\psi_{2}) \rbrace \\
&= \mathscr{F}_{S}( \CfdRdxR ).
\end{aligned} \right.
\end{equation}
Members of this class act as members of $\VRd$ spatially and as members of $ \CfdR$ temporally. That is, if $\psi \in \VRdCfdR$,  then for every $x \in \Rd$, $\psi(x , \cdot) \in C_{FD}(\R)$ and for every $t \in \R$, $\psi( \cdot , t ) \in \VRd$. It is immediate that this set of functions is a complex vector space.  Every function of the form $\varphi \boxtimes \theta $, with $\varphi \in \mathscr{V}(\Rd)$ and $\theta \in \CfdR $ is a member of this class, as well as any finite linear combination of functions of this form. The notation $\VRdCfdR $ has a deep inspiration in the theory of Nuclear spaces: the notation $E\hat{\otimes}F$, when $E$ and $F$ are general topological vector spaces, is used to represent a completition, under suitable topologies, of the space $E\otimes F$ of finite linear combinations of tensor products (see \citet[Part III]{treves1967topological} or \citet{grothendieck1955produits}). We will not enter in those details explicitly, and we will simply work with definition \eqref{Eq:DefVRdCfdR} and its notation, which would be fully justified if it turns out that the space $\VRd$ is nuclear \citep[Proposition 50.7]{treves1967topological}.

We endow $\VRdCfdR$ with the topology induced by the directed family of semi-norms:
\begin{equation}
\label{Eq:SeminormsVRdCfdR}
p_{n_{S} , n_{T}}(\psi) = \sup_{(\xi,t) \in \RdxR}\left| (1+|\xi|^{2})^{n_{S}}(1+t^{2})^{n_{T}}\mathscr{F}^{-1}_{S}(\psi)(\xi,t) \right|, \quad n_{S}, n_{T} \in \mathbb{N}. 
\end{equation}
The spatial Fourier transform (and its inverse) defines a continuous linear functional from $C_{FD}(\RdxR)$ to $\VRdCfdR$ (or from $\VRdCfdR$ to $\CfdRdxR$, as pleasure). From the continuity of $\mathscr{F}_{S}$ and the completeness of $C_{FD}(\RdxR)$, we conclude that $\VRdCfdR$ is a Fréchet space. The dual space of $\VRdCfdR$ is denoted by $\VRdCfdRDual$. A linear functional $T : \VRdCfdR \to \mathbb{C}$ is a member of $\VRdCfdRDual$ if and only if there exist $C > 0 $ and $N_{S},N_{T} \in \mathbb{N}$ such that
\begin{equation}
\label{Eq:CriterionTVRdCfdRDual}
| \langle T , \psi \rangle | \leq C \sup_{(\xi,t) \in \RdxR}\left|(1+|\xi|^{2})^{N_{S}}(1+t^{2})^{N_{T}}\mathscr{F}_{S}^{-1}(\psi)(\xi,t)  \right|, \quad \forall \psi \in \VRdCfdR.
\end{equation}
The spatial Fourier transform $\mathscr{F}_{S}$ over $\VRdCfdRDual$ is defined as the adjoint of the spatial Fourier transform over $C_{FD}(\RdxR)$, whose range is the space $\VRdCfdR$. We obtain thus a continuous linear operator $\mathscr{F}_{S} :  \VRdCfdRDual \to \mathscr{M}_{SG}(\RdxR)$. It is then concluded that $\VRdCfdRDual = \mathscr{F}_{S}\left(\MsgRdxR \right)$, the space of tempered distributions whose spatial Fourier transforms are slow-growing measures over $\RdxR$.

We consider now spaces of test-functions which are conceived to work over the positive-time. We define the next space of test-functions:
\begin{equation}
\label{Eq:DefCfdRdxRp}
C_{FD}(\RdxRp) := \lbrace \psi \in C(\RdxRp) \ \big| \sup_{(\xi,t)\in \RdxRp} \left| (1+|x|^{2})^{n_{S}}(1+t^{2})^{n_{T}} \psi(\xi,t)  \right| < \infty, \forall n_{S},n_{T} \in \mathbb{N} \rbrace.
\end{equation}
This space can be equivalently defined as the space of restrictions of functions in $C_{FD}(\RdxR)$ to the subset $\RdxRp$. $C_{FD}(\RdxRp)$ is endowed with the topology induced by the family of semi-norms:
\begin{equation}
\label{Eq:SeminormsCfdRdxRp}
p_{n_{S},n_{T}}(\psi) = \sup_{(\xi,t) \in \RdxRp}\left|(1+|\xi|^{2})^{n_{S}}(1+t^{2})^{n_{T}}\psi(\xi,t) \right| , \quad n_{S}, n_{T} \in \mathbb{N}.
\end{equation}
$C_{FD}(\RdxRp)$ is a Fréchet space. We denote by $C_{FD}'(\RdxRp)$ its dual space. A linear functional $T : C_{FD}(\RdxRp) \to \mathbb{C}$ is in $\CfdRdRpDual$  if and only if there exist $C > 0 $ and $N_{S},N_{T} \in \mathbb{N} $ such that
\begin{equation}
\label{Eq:CriterionTinCfdRpDual}
\left| \langle T , \psi \rangle \right| \leq C \displaystyle\sup_{(\xi,t) \in \RdxRp}\left|(1+|\xi|^{2})^{N_{S}}(1+t^{2})^{N_{T}}\psi(\xi,t)  \right|, \quad \forall \psi \in C_{FD}(\RdxRp).
\end{equation}

Let us denote $\MsgRdxRp$  the space of slow-growing complex measures over $\RdxR$ with support contained in $\RdxRp$. We could have defined $\mathscr{M}_{SG}(\RdxRp)$ as a space of measures over $\RdxRp$ without concerning on what happens over the \textit{negative-time}, but it is actually easier to work with measures defined over the whole space $\RdxR$ but for which we only look at their behaviours over the subset $\RdxRp$. The next Proposition follows from Riesz Representation Theorem \ref{Theo:RieszReperesentationDualCfdMsg}.
\begin{Prop}
\label{Prop:MsgRdRpDualCfdRdRp}
$\mathscr{M}_{SG}(\RdxRp) = C_{FD}'(\RdxRp)$.
\end{Prop}
\textbf{Proof: } The inclusion $\mathscr{M}_{SG}(\RdxRp) \subset C_{FD}'(\RdxRp)$ is immediate. Let $T \in C_{FD}'(\RdxRp)$. We extend the domain of definition of $T$ so it can be applied to every function in $\CfdRdxR$, through
\begin{equation}
\label{Eq:DefTRestrictionMeasure}
\langle T , \psi \rangle := \langle T , \psi\big|_{\RdxRp} \rangle, \quad \forall \psi \in \CfdRdxR.
\end{equation}
Since $T \in \CfdRdRpDual$, there exist $C >0 $ and $N_{S},N_{T} \in \mathbb{N}$ such that
\begin{equation}
\label{Eq:BoundingArgumentTRiezsMeasure}
\left.\begin{aligned}
\left| \langle T , \psi \rangle \right| &=  \left| \langle T , \psi\big|_{\RdxRp} \rangle \right| \\
 &\leq C \displaystyle\sup_{(\xi,t) \in \RdxRp}\left|(1+|\xi|^{2})^{N_{S}}(1+t^{2})^{N_{T}}\psi(\xi,t)  \right| \\
 &\leq C \displaystyle\sup_{(\xi,t) \in \RdxR}\left|(1+|\xi|^{2})^{N_{S}}(1+t^{2})^{N_{T}}\psi(\xi,t)  \right|.
\end{aligned}\right.
\end{equation}
This proves that the extension of $T$ to the space $C_{FD}(\RdxR)$ is in  $\CfdRdxRDual$ and hence by Riesz Representation Theorem \ref{Theo:RieszReperesentationDualCfdMsg} there exists a slow-growing measure $\mu \in \mathscr{M}_{SG}(\RdxR)$ such that
\begin{equation}
\label{Eq:TRepresentedMeasureSupportRdRp}
\langle T , \psi\big|_{\RdxRp} \rangle =  \int_{\RdxR} \psi(x,t) d\mu(x,t),\quad \forall \psi \in C_{FD}(\RdxR).
\end{equation}    
It can be concluded that the support of $\mu$ is contained in $\RdxRp$ by considering that $\langle T , \psi \rangle = 0 $ for every $\psi \in \CfdRdxR$ such that its support is contained in $\RdxRneg$. This proves that $C_{FD}'(\RdxRp) = \mathscr{M}_{SG}(\RdxRp)$. $ \blacksquare $

We consider now the following space of test-functions:
\begin{equation}
\label{Eq:DefVRdCfdRp}
\left. \begin{aligned}
\VRdCfdRp  &:= \lbrace \psi \in C(\RdxRp) \ \big| \ \exists \psi_{2} \in C_{FD}(\RdxRp) \hbox{ such that } \psi = \mathscr{F}_{S}(\psi_{2}) \rbrace \\
&= \mathscr{F}_{S}( \CfdRdxRp ).
\end{aligned} \right.
\end{equation}
Members of this class satisfy completely analogous properties to those of the class $\VRdCfdR$, with the only difference residing in the domain of definition of the functions. The notation of \eqref{Eq:DefVRdCfdRp} is also inspired by the theory of Nuclear spaces, although not justified.  The space $\VRdCfdRp$ is also endowed with an analogous topology  as the space $\VRdCfdR$, using the supremum over $\RdxRp$ rather than over $\RdxR$ in the definition of the semi-norms \eqref{Eq:SeminormsVRdCfdR}. $\VRdCfdRp$ is a Fréchet space. The spatial Fourier transform defines a continuous linear operator from $\VRdCfdRp$ to $\CfdRdxRp$ and vice-versa.  The dual space of $\VRdCfdRp$ is noted $\VRdCfdRpDual$. A linear functional $T : \VRdCfdRp \to \mathbb{C}$ is a member of $\VRdCfdRpDual$ if and only if there exist $C > 0 $ and $N_{S},N_{T} \in \mathbb{N}$ such that
\begin{equation}
\label{Eq:CriterionTVRdCfdRpDual}
| \langle T , \psi \rangle | \leq C \sup_{(\xi,t) \in \RdxRp}\left|(1+|\xi|^{2})^{N_{S}}(1+t^{2})^{N_{T}}\mathscr{F}_{S}^{-1}(\psi)(\xi,t)  \right|, \quad \forall \psi \in \VRdCfdRp.
\end{equation}
We conclude that $\VRdCfdRpDual = \mathscr{F}_{S}(\MsgRdxRp)$. Since $\MsgRdxRp \subset \MsgRdxR$, we also conclude that $\VRdCfdRpDual$ is a subspace of $ \VRdCfdRDual$.

We remark that both spaces $\VRdCfdRDual$ and $\MsgRdxR$ are subspaces of $\TemperedRdxR$. It turns out that every distribution in $\mathscr{M}_{SG}(\RdxR)$ or in $\VRdCfdRpDual$ can be differentiated any number of times, and that the spatial, temporal and spatio-temporal Fourier transforms can be applied freely.

\subsection{Operators defined through a symbol}
\label{Sec:OperatorsSymbol}

Here we work with a generic $d \in \mathbb{N}_{*}$, not necessarily in a spatio-temporal framework. Let $g : \Rd \to \mathbb{C}$. Let us denote $g_{R}$ its real part and $g_{I}$ its imaginary part. We say that $g$ is a \textit{symbol function} if it is measurable, polynomially bounded and Hermitian (that is, $g_{R}$ is even and $g_{I}$ is odd). If $g$ is a symbol function, we define the linear operator $\mathcal{L}_{g} : \VRdDual \mapsto \VRdDual$ given by
\begin{equation}
\label{Eq:DefLg}
\mathcal{L}_{g}(T) = \mathscr{F}^{-1}( g \mathscr{F}(T) ), \quad \forall T \in \VRdDual. 
\end{equation}
This operator is well-defined since the multiplication of a measurable and polynomially bounded function with a slow-growing measure is a slow-growing measure. In addition, the Hermitianity of $g$ guarantees that $\mathcal{L}_{g}$ is a real operator. Every differential operators with constant coefficients is of this form, obtained when $g$ is an Hermitian polynomial. Other pseudo-differential operators such as $(\kappa^{2}-\Delta)^{\frac{\alpha}{2}}$ for $\alpha \in \R $ and $\kappa > 0 $ (or $\alpha \geq 0 $ and $\kappa = 0 $) can be also obtained through this method through the function $g(\xi) = (\kappa^{2} + |\xi|^{2})^{\frac{\alpha}{2}} $. We remark that when $g$ is continuous, $\mathcal{L}_{g}$ can be identified as the adjoint of the operator $\mathscr{F}( g\mathscr{F}^{-1}( \cdot ) )$, for which it is quite immediate to see that it is a continuous linear operator from $\VRd$ to $\VRd$.

We remark that in the case where $|g| \geq \frac{1}{p}$ for a strictly positive polynomial $p : \Rd \to \Rp_{*}$, the operator $\mathcal{L}_{g}$ is bijective, with $\mathcal{L}_{g}^{-1} = \mathcal{L}_{\frac{1}{g}}$. In such a case, if we consider a PDE of the form 
\begin{equation}
\label{Eq:PDELgU=X}
\mathcal{L}_{g}U = X, \quad X \in \VRdDual,
\end{equation}
we conclude that there exists a unique solution in $\VRdDual$ given by $U = \mathcal{L}_{\frac{1}{g}}(X)$. A more general treatment for equations of the form \eqref{Eq:PDELgU=X} can be found in \citet[Section 4.6]{carrizov2019development}.

In a spatio-temporal setting, for a spatial symbol function $g$, the operator $\mathcal{L}_{g}$ will denote the operator $\mathscr{F}_{S}^{-1}(g\mathscr{F}_{S}(\cdot))$, which is defined over $\VRdCfdRDual$. The PDE \eqref{Eq:ThePDE} can thus be seen as an equation with an operator defined through a spatio-temporal symbol, $\frac{\partial}{\partial t} + \mathcal{L}_{g} = \mathscr{F}^{-1}\left(   (i\omega + g )\mathscr{F}(\cdot) \right)$.

\subsection{Càdlàg-in-time representations}
\label{Sec:CadlagInTimeRepresentations}

In this section we make precise what do we mean with distributions having a \textit{càdlàg-in-time behaviour}. This definition will be restrained in this work to the cases of the spaces of distributions presented in Section \ref{Sec:ConvinientTestFunctions}, although a more general treatment can be done.

\begin{Def}
\label{Def:TCadlagInTime}
Let $T \in \VRdCfdRDual$ (resp. in $\MsgRdxR$). We say that $T$ has a càdlàg behaviour over $\R$, or that it has a càdlàg-in-time representation, if there exists a family of spatial distributions $(T_{t})_{t \in \R} \subset \VRdDual$ (resp.  $(T_{t})_{t \in \R} \subset \MsgRd$) such that
\begin{itemize}
\item for all $\varphi \in \VRd$ (resp. in $\CfdRd$) the function $t \mapsto \langle T_{t} , \varphi \rangle$ is a càdlàg function defining a slow-growing measure over $\R$,
\item for all $\varphi \in \VRd$ (resp. in $\CfdRd$) and for all $\theta \in \CfdR$, it holds that
\begin{equation}
\label{Eq:TcadlagInTimeIntegral}
\langle T  , \varphi \boxtimes \theta \rangle = \int_{\R}\langle T_{t} , \varphi \rangle \theta(t) dt.
\end{equation}
\end{itemize}
The family of distributions $(T_{t})_{t \in \R}$ is called the càdlàg-in-time representation of $T$.
\end{Def}

One can easily conclude that if $T$ is a distribution in $\VRdCfdRDual$ (resp. in $\MsgRdxR$) having a càdlàg-in-time representation, its representation must be unique. Indeed, if $(T_{t})_{t \in \R}$ and $(\tilde{T}_{t})_{t \in \R}$ are two càdlàg-in-time representations of $T$, then for every $\varphi$ in $\VRdCfdR$ (resp. in $\CfdRdxR$),  $ \int_{\R} \langle T_{t} - \tilde{T}_{t} , \varphi \rangle \theta(t) dt = 0 $ for every $\theta \in \CfdR$, from which it can be concluded that $t \mapsto \langle T_{t} - \tilde{T}_{t} , \varphi \rangle$ equals $0$ almost everywhere. Since this last function is càdlàg, it must be null everywhere. We remark that in the case where $T \in \VRdCfdRpDual$ (resp. in $\MsgRdxRp$) has a càdlàg-in-time representation, then $T_{t} = 0 $ for every $t \in \Rneg$.

Let us introduce the next auxiliary definition.

\begin{Def}
\label{Def:ConvergingDeltaRight}
Let $t_{0} \in \R$. Consider a sequence of positive functions $(\theta_{n}^{(t_{0})})_{n \in \mathbb{N}} \subset \DR$ such that $\int_{\R}\theta_{n}^{(t_{0})}(t)dt =1 $ for all $n \in \mathbb{N}$ and such that there exists a sequence of strictly positive real numbers $(a_{n})_{n \in \mathbb{N}} \subset \Rp_{*}$  which decreases to $0$ and such that $\supp( \theta_{n}^{(t_{0})} ) \subset \left[ t_{0} , t_{0} + a_{n} \right] $ for all $n \in \mathbb{N}$. We say then that $(\theta_{n}^{(t_{0})})_{n \in \mathbb{N}}$ converges to $\delta_{t_{0}}$ from the right, denoted $\theta_{n}^{(t_{0})} \to \delta_{t_{0}}^{+}$. The definition of a sequence converging to $\delta_{t_{0}}$ from the left, denoted $\theta_{n}^{(t_{0})} \to \delta_{t_{0}}^{-}$, is analogous.
\end{Def}

It is clear that if $ f : \R \to \mathbb{C}$ is a càdlàg function then for all $t \in \R$, $f(t) = \lim_{n \to \infty} \int_{\R}f(s)\theta_{n}^{(t)}(s)ds$ for any sequence $(\theta_{n}^{(t)})_{n \in \mathbb{N}}$ converging to $\delta_{t}$ from the right. Hence, if $T$ is a distribution in any of the spaces considered in Definition \ref{Def:TCadlagInTime} and being càdlàg-in-time, the càdlàg-in-time representation is determined by
\begin{equation}
\label{Eq:TtAsLimit}
\langle T_{t} , \varphi \rangle = \lim_{n \to \infty} \langle T , \varphi \boxtimes \theta_{n}^{(t)} \rangle, \quad \forall t \in \R,
\end{equation}
for any sequence $(\theta_{n}^{(t)})_{n \in \mathbb{N}}$ converging to $\delta_{t}$ from the right, and for any spatial test-function $\varphi$ in the corresponding space. The next Proposition is quite obvious.

\begin{Prop}
\label{Prop:CadlagRepresentationInterchangeSpatialFourier}
Let $T \in \MsgRdxR$. Then, $T$ has a càdlàg-in-time representation if and only if $\mathscr{F}_{S}(T) \in \VRdCfdRDual$ has a càdlàg-in-time representation.
\end{Prop}

\textbf{Proof: } It suffices to take $(\mathscr{F}_{S}(T_{t}))_{t  \in \R}$ as the family defining the càdlàg-in-time representation of $\mathscr{F}_{S}(T)$. $\blacksquare$

Let $T$ be in $\VRdCfdRDual$ or in $\MsgRdxR$. One crucial fact about $T$ is that it acts as a \textit{slow-growing measure in time}. This fact is evident for the case $T \in \MsgRdxR$. For $T \in \VRdCfdRDual$, what we mean with $T$ being a slow-growing measure in time is that for every spatial test-function $\varphi \in \VRd$, the application $\theta \in \CfdR \mapsto \langle T , \varphi \boxtimes \theta \rangle$  defines a continuous linear functional over $\CfdR$ and hence it is a slow-growing measure over $\R$. This can be concluded immediately from criteria \eqref{Eq:CriterionTVRdCfdRDual}. The fact that $T$ acts temporally as a measure implies that we can construct integrals with respect to its time component. To be precise, for every spatial test-function $\varphi$ in a corresponding space, we can extend the domain of definition of the application $\theta \mapsto \langle T , \varphi \boxtimes \theta \rangle  $ to every measurable and bounded function with fast decreasing behaviour over $\R$. In particular, for any $A \in \BoundedBorelRd$, we have the right to write:
\begin{equation}
\label{Eq:TvarphiBoxtimes1A}
\langle T , \varphi\boxtimes \mathbf{1}_{A} \rangle.
\end{equation}

The next Theorem generalises in some sense the condition we can testify in dimension $1$: that the temporal \textit{primitive} of a measure in $\MR$ can be identified with a càdlàg function (Eq. \eqref{Eq:CadlagPrimitiveMeasure1D}).

\begin{Theo}
\label{Theo:CadlagRepresentation}
Let $\mathscr{U}$ denoting  any of the spaces $\VRdCfdRDual$, $\VRdCfdRpDual$, $\MsgRdxR$ or $\MsgRdxRp$. Let $T \in \mathscr{U}$ be such that $\frac{\partial T}{\partial t} \in \mathscr{U}$. Then, $T$ has a càdlàg-in-time representation.
\end{Theo}

\textbf{Proof:} We will only prove the case where $T$ and $\frac{\partial T}{\partial t}$ are in $\MsgRdxR$. The other cases follow immediately form $\MsgRdxRp \subset \MsgRdxR$ and from Proposition \ref{Prop:CadlagRepresentationInterchangeSpatialFourier}. We will denote for simplicity $\mu = \frac{\partial T}{\partial t} \in \MsgRdxR$. In this proof we will extensively use that fact that if $\theta_{1} , \theta_{2} \in \DR$ are two temporal functions such that $\int_{\R} \theta_{1}(t)dt = \int_{\R} \theta_{2}(t)dt = 1$, then the function $\theta_{1} - \theta_{2}$ has a unique primitive in $\DR$, given by $t \in \R \mapsto \int_{-\infty}^{t}\theta_{1}(u) - \theta_{2}(u) du$. In addition, this primitive has its support contained in the support of $\theta_{1} - \theta_{2}$. We will denote this primitive $\int (\theta_{1} - \theta_{2} )$.  Hence, we will always have for functions of this form,
\begin{equation}
\label{Eq:TappliedDifferenceTestFunctions}
\langle T , \varphi \boxtimes (\theta_{1}- \theta_{2})  \rangle = - \langle \frac{\partial T}{\partial t} , \varphi \boxtimes \int( \theta_{1}-\theta_{2} ) \rangle = -\int_{\Rd\times\supp(\theta_{1} - \theta_{2} )} \varphi(\xi)\int_{-\infty}^{t}\theta_{1}(u) - \theta_{2}(u) du d\mu(\xi,t).
\end{equation}

Let $\varphi \in \CfdRd$. Let $t \in \R$ and let $(\theta_{n}^{(t)})_{n \in \mathbb{N}} \subset \DR$ be a sequence convergning to $\delta_{t}$ from the right, for which we will suppose for simplicity that $\supp(\theta_{n}^{(t)}) \subset \left[ t , t + \frac{1}{n} \right]$ for all $n \in \mathbb{N}_{*}$.  We consider thus the sequence of complex numbers $ ( \langle T , \varphi \boxtimes \theta_{n}^{(t)} \rangle)_{n \in \mathbb{N}} $. For $n,m \in \mathbb{N}_{*}$, the function $\theta_{n}^{(t)} - \theta_{m}^{(t)}$ has null total integral, and hence we can apply principle \eqref{Eq:TappliedDifferenceTestFunctions} to it. We have thus,
\begin{equation}
\label{Eq:TVarphiRestaThetasConvergesZero}
\left. \begin{aligned}
\left| \langle T , \varphi \boxtimes \theta_{n}^{(t)} \rangle - \langle T , \varphi \boxtimes \theta_{m}^{(t)} \rangle \right| &= \left| \langle T , \varphi \boxtimes \left(\theta_{n}^{(t)} - \theta_{m}^{(t)} \right) \rangle \right|.  \\
&=  \left| \int_{\Rd\times\left(t , t + \frac{1}{n} \vee \frac{1}{m} \right) } \varphi(\xi) \int_{-\infty}^{s}  \theta_{n}^{(t)}(u) - \theta_{m}^{(t)}(u) du d\mu(\xi,s) \right| \\
&\leq   \int_{\Rd\times\left(t , t + \frac{1}{n} \vee \frac{1}{m} \right) } 2|\varphi(\xi)| d|\mu|(\xi,s) \to 0 \quad \hbox{ as } n,m \to \infty. 
\end{aligned} \right.
\end{equation}
Here we have used that $\left|\int_{-\infty}^{s}\theta_{n}^{(t)}(u) - \theta_{m}^{(t)}(u) du \right| \leq  \int_{\R}\theta_{n}^{(t)}(u) + \theta_{m}^{(t)}(u) du  = 2$. The convergence to $0$ is justified since the set $\Rd\times\left(t , t + \frac{1}{n} \vee \frac{1}{m} \right)$ decreases to $\emptyset$ as $n,m \to \infty$ and since $\varphi \in \CfdRd$. This proves that the sequence $ ( \langle T , \varphi \boxtimes \theta_{n}^{(t)} \rangle)_{n \in \mathbb{N}} $ is a Cauchy sequence and hence it converges to a limit which we will denote $\langle T_{t} , \varphi \rangle$. An argument following the same procedures as in \eqref{Eq:TVarphiRestaThetasConvergesZero} can be used to prove that the limit does not depend on the sequence $(\theta_{n}^{(t)})_{n \in \mathbb{N}}$ converging to $\delta_{t}$ chosen (replace $\theta_{m}^{(t)}$ with another sequence, the arguments still hold).

Since $T_{t}$ satisfies \eqref{Eq:TtAsLimit} for every $t \in \R$, it is immediate that the application $\varphi \in \CfdRd \mapsto \langle T_{t} , \varphi \rangle$ is linear. Consider a sequence $(\theta_{n}^{(t)})_{n \in \mathbb{N}}$ converging to $\delta_{t}$ from the right side, and let us consider a positive function $\vartheta_{t} \in \DR$ such that $\int_{\R}\vartheta_{t}(u)du=1$, $ \supp( \vartheta_{t}) \subset \left[ t , t+1 \right]$, and such that $\| \vartheta \|_{\infty} := \sup_{u \in \R}| \vartheta_{t}(u) |$ does not depend on $t$ (we can always find such a function). One has then
\begin{equation}
\label{Eq:BoundingTtVarphiThetant}
\left. \begin{aligned}
\left| \langle T , \varphi \boxtimes \theta_{n}^{(t)} \rangle \right| &\leq | \langle T , \varphi \boxtimes (\theta_{n}^{(t)} - \vartheta_{t}) \rangle | + | \langle T , \varphi \boxtimes \vartheta_{t} \rangle | \\
&=  \left| \int_{\Rd\times\left[ t , t+ 1 \right]} \varphi(\xi)\int_{-\infty}^{s} \theta_{n}^{(t)}(u) - \vartheta_{t}(u) du d\mu(\xi,s) \right| + \left| \int_{\Rd \times\left[t , t+1 \right]}\varphi(\xi)\vartheta_{t}(s) dT(\xi,s) \right|.
\end{aligned} \right.
\end{equation}
Let $N^{\mu}_{S} , N^{\mu}_{T}, N_{S}^{T}, N_{T}^{T} \in \mathbb{N}$ be such that the measures $ (1+|\xi|^{2})^{-N^{\mu}_{S}}(1 + t^{2} )^{-N^{\mu}_{T}}\mu $ and $(1+|\xi|^{2})^{-N^{\mu}_{S}}(1 + t^{2} )^{-N^{\mu}_{T}}T$ are finite. Using Lemma \ref{Lemma:ConvexityInequalities} we conclude the inequalities
\begin{equation}
\left. \begin{aligned}
\Big| \int_{\Rd\times\left[ t , t+ 1 \right]} \varphi(\xi)&\int_{-\infty}^{s} \theta_{n}^{(t)}(u) - \vartheta_{t}(u) du d\mu(\xi,s) \Big| \leq \int_{\Rd\times\left[t,t+1\right]}|\varphi(\xi)|2 (1+s^{2})^{N^{\mu}_{T}} \frac{d|\mu|(\xi,s)}{(1+s^{2})^{N^{\mu}_{T}}} \\
&\leq 2^{N_{T}^{\mu}}\sup_{\xi \in \Rd}\left| (1+|\xi|^{2})^{N_{S}^{\mu}} \varphi(\xi) \right| \int_{\RdxR} \frac{d|\mu|(\xi,s)}{(1+|\xi|^{2})^{N_{S}^{\mu}}(1+s^{2})^{N_{T}^{\mu}}} \left( 3^{N_{T}^{\mu}} + 2^{N_{T}^{\mu}} t^{2 N_{T}^{\mu}} \right),
\end{aligned} \right.
\end{equation}
and
\small
\begin{equation}
\left. \begin{aligned}
 \Big| \int_{\Rd\times\left[ t , t+1 \right]}\varphi(\xi)&\vartheta_{t}(s) dT(\xi,s)  \Big| \leq  \int_{\Rd\times\left[t , t+1 \right]} |\varphi(\xi)| |\vartheta_{t}(s)|d|T|(\xi,s) \\
&\leq \| \vartheta \|_{\infty}\int_{\Rd\times\left[t,t+1\right]} |\varphi(\xi)|(1+|\xi|^{2})^{N_{S}^{T}} (1+s^{2})^{N_{T}^{T}} \frac{d|T|(\xi,s)}{(1+|\xi|^{2})^{N_{S}^{T}} (1+s^{2})^{N_{T}^{T}}} \\
&\leq  2^{N_{T}^{T}-1}\| \vartheta \|_{\infty} \sup_{\xi \in \Rd}\left| (1+|\xi|^{2})^{N_{S}^{T}}\varphi(\xi) \right| \int_{\RdxR}\frac{d|T|(\xi,s)}{(1+|\xi|^{2})^{N_{S}^{T}} (1+s^{2})^{N_{T}^{T}}}\left( 3^{N_{T}^{T}} + 2^{N_{T}^{T}}t^{2N_{T}^{T}} \right). 
\end{aligned} \right.
\end{equation}
\normalsize
We conclude thus
\small
\begin{equation}
\label{Eq:TtIsCondinuousLinearFunctional}
\left. \begin{aligned}
\left| \langle T_{t} , \varphi \rangle \right| &= \lim_{n \to \infty} \left| \langle T , \varphi \boxtimes \theta_{n}^{(t)} \rangle \right| \\
&\leq  \Bigg[ 2^{N_{T}^{\mu}} \int_{\RdxR} \frac{d|\mu|(\xi,s)}{(1+|\xi|^{2})^{N_{S}^{\mu}}(1+s^{2})^{N_{T}^{\mu}}} \left( 3^{N_{T}^{\mu}} + 2^{N_{T}^{\mu}} t^{2 N_{T}^{\mu}} \right) \\
&\quad + 2^{N_{T}^{T}-1}\| \vartheta \|_{\infty} \int_{\RdxR}\frac{d|T|(\xi,s)}{(1+|\xi|^{2})^{N_{S}^{T}}(1+s^{2})^{N_{T}^{T}}}\left( 3^{N_{T}^{T}} + 2^{N_{T}^{T}}t^{2N_{T}^{T}} \right) \Bigg]\sup_{\xi \in \Rd}\left| (1+|\xi|^{2})^{N_{S}^{\mu}\vee N_{T}^{T}} \varphi(\xi) \right|.
\end{aligned} \right.
\end{equation}
\normalsize
This inequality proves two things: first, that the function $t \in \R \mapsto \langle T_{t} , \varphi \rangle $ is polynomially bounded for every $\varphi \in \CfdRd$, and second that for every $t \in \R$ the application $\varphi \in \CfdRd \mapsto \langle T_{t} , \varphi \rangle$ defines a continuous linear functional over $\CfdRd$, hence it defines a slow-growing measure over $\Rd$. Let us prove that the function $t \in \R \mapsto \langle T_{t} , \varphi \rangle$ is càdlàg. Let $t,s \in \R$ such that $t > s$. Let $(\theta_{n}^{(s)})_{n \in \mathbb{N}}$ converging to $\delta_{s}$ from the right and let $(\theta_{n}^{(t)})_{n \in \mathbb{N}}$ converging to $\delta_{t}$ from the right. Let us make explicit the expression $\langle T_{t} - T_{s} , \varphi \rangle$. We have
\begin{equation}
\left. \begin{aligned}
\langle T_{t} , \varphi \rangle - \langle T_{s} , \varphi \rangle &= \langle T_{t} - T_{s} , \varphi \rangle \\
&= \lim_{n \to \infty} \langle T , \varphi \boxtimes (\theta_{n}^{(t)} - \theta_{n}^{(s)}) \rangle \\
&= \lim_{n \to \infty} - \int_{\RdxR} \varphi(\xi) \int_{-\infty}^{u} \theta_{n}^{(t)}(v) - \theta_{n}^{(s)}(v) dv d\mu(\xi,u) \\
&=  \int_{\RdxR} \varphi(\xi)\mathbf{1}_{\left( s , t \right]}(u) d\mu(\xi, u).
\end{aligned} \right.
\end{equation}
The last limit is justified using Dominated Convergence Theorem. Indeed, one has $$ \lim_{n \to \infty}\varphi(\xi) \int_{-\infty}^{u} \theta_{n}^{(t)}(v) - \theta_{n}^{(s)}(v) dv \to -\varphi(\xi)\mathbf{1}_{\left( s , t \right]}(u), \quad \forall (\xi,u) \in \RdxR,$$ and this point-wise convergence is dominated by the function $(\xi,u) \mapsto 2|\varphi(\xi)|\mathbf{1}_{\left[s , t + 1\right)}(u) $ which is integrable with respect to $|\mu|$. We write then
\begin{equation}
\label{Eq:Tt-Ts}
\langle T_{t} - T_{s} , \varphi \rangle = \langle \mu , \varphi \boxtimes \mathbf{1}_{\left( s , t \right]} \rangle,
\end{equation}
from which it is clear that the function $t \mapsto \langle T_{t} , \varphi \rangle$ is càdlàg. Since this càdlàg function is polynomially bounded, it defines a slow-growing measure over $\R$, with the integral $\int_{\R} \langle T_{t} , \varphi \rangle \theta(t)dt$ being well-defined for every $\theta \in \CfdR$.

We finish by proving that $\int_{\R} \langle T_{t} , \varphi \rangle\theta(t)dt = \langle T , \varphi \boxtimes \theta \rangle$ for all $\varphi \in \CfdRd$ and $\theta \in \CfdR$. For that, we will consider for every $t \in \R$ a particular sequence of functions converging to $\delta_{t}$ from the right given by
\begin{equation}
\label{Eq:SeqThetas}
\theta_{n}^{(t)}(s) = a_{n}e^{-\frac{1}{1-\left| 2n\left( s-\frac{1}{2} - \left(t + \frac{1}{2n} \right) \right) \right|^{2}}}\mathbf{1}_{\left( t , t + \frac{1}{n}  \right) }(s), \quad n \in \mathbb{N}_{*},
\end{equation}
where $a_{n} > 0 $ is a normalising number. This sequence has the particularity that if we regard the function $t \mapsto \theta_{n}^{(t)}(s)$ for a fixed $s \in \R$, the obtained sequence of functions approaches $\delta_{s}$ from the left. We also remark that the function $(t,s) \in \R\times\R \mapsto \theta_{n}^{(t)}(s)$ is measurable. From  \eqref{Eq:TtAsLimit}, one has that
\begin{equation}
\label{Eq:TtThetaLimit}
\langle T_{t} , \varphi \rangle\theta(t) = \lim_{n \to \infty} \langle T , \varphi \boxtimes \theta_{n}^{(t)} \rangle \theta(t), \quad \forall t \in \R.
\end{equation}
In addition, following \eqref{Eq:TtIsCondinuousLinearFunctional} one has that the sequence of functions $( t \mapsto \langle T , \varphi \boxtimes \theta_{n}^{(t)} \rangle )_{n \in \mathbb{N}_{*}}$ is uniformly bounded by a polynomial, hence integrable when multiplied by $\theta$. The point-wise convergence  \eqref{Eq:TtThetaLimit} is thus dominated by an integrable function. We have thus
\begin{equation}
\int_{\R}\langle T , \varphi \rangle \theta(t) dt = \int_{\R} \lim_{n \to \infty} \langle T , \varphi \boxtimes \theta_{n}^{(t)} \rangle \theta(t) dt = \lim_{n \to \infty} \int_{\R} \langle T , \varphi \boxtimes \theta_{n}^{(t)} \rangle \theta(t) dt.
\end{equation}
Using again Lemma \ref{Lemma:ConvexityInequalities}, we consider that for all $(\xi,t,s) \in \RdxR\times\R$, it holds that
\begin{equation}
\left. \begin{aligned}
|\varphi(\xi)||\theta_{n}^{(t)}(s)| |\theta(t) | &= |\varphi(\xi)||\theta_{n}^{(t)}(s)| |\theta(t) | \mathbf{1}_{\left[ t , t+1 \right]}(s)\frac{(1+s^{2})^{N_{T}^{T}}}{(1+s^{2})^{N_{T}^{T}}} \\
&\leq a_{n}2^{N_{T}^{T}-1}|\varphi(\xi)|(3^{N_{T}^{T}} + 2^{N_{T}}t^{2N_{T}} )|\theta(t)|\frac{1}{(1+s^{2})^{N_{T}^{T}}},
\end{aligned} \right.
\end{equation}  
and the last function is integrable with respect to the measure $d|T|(\xi,s)dt$ over $\RdxR\times\R$. By Fubini's Theorem, it follows that
\begin{equation}
\int_{\R} \langle T , \varphi \boxtimes \theta_{n}^{(t)} \rangle \theta(t) dt = \int_{\R} \int_{\RdxR} \varphi(\xi) \theta_{n}^{(t)}(s) \theta(t) dT(\xi,s)dt  =  \int_{\RdxR} \varphi(\xi) \int_{\R}  \theta_{n}^{(t)}(s) \theta(t) dt dT(\xi,s).
\end{equation}
Since for any fixed $s \in \R$, $(t \mapsto \theta_{n}^{(t)}(s) ) \to \delta_{s}^{-}$ and $\theta$ is continuous, it follows that $\lim_{n \to \infty} \int_{\R} \theta(t) \theta_{n}^{(t)}(s) dt = \theta(s) $. Using again lemma \ref{Lemma:ConvexityInequalities}, one has for all $(\xi,s) \in \RdxR$,
\begin{equation}
\left. \begin{aligned}
\left|\varphi(\xi)\int_{\R}\theta(t)\theta_{n}^{(t)}(s)dt\right| &= \left|\varphi(\xi)\frac{(1+s^{2})^{N_{T}^{T}}}{(1+s^{2})^{N_{T}^{T}}}\int_{\left[ s-1,s \right]}\theta(t)\theta_{n}^{(t)}(s)dt\right| \\
&\leq 2^{N_{T}^{T}-1}\frac{|\varphi(\xi)|}{(1+s^{2})^{N_{T}^{T}}}  \int_{\left[ s-1 , s \right]} \left(  3^{N_{T}^{T}} + 2^{N_{T}^{T}}t^{2N_{T}^{T}} \right)|\theta(t) | \theta_{n}^{(t)}(s)dt \\
&\leq 2^{N_{T}^{T}-1} \sup_{t \in \R}\left|  \left(  3^{N_{T}^{T}} + 2^{N_{T}^{T}}t^{2N_{T}^{T}} \right)\theta(t) \right|\frac{|\varphi(\xi)|}{(1+s^{2})^{N_{T}^{T}}}.
\end{aligned} \right.
\end{equation}
This last function is integrable with respect to $|T|$ over $\RdxR$. We conclude
\begin{equation}
\lim_{n \to \infty} \int_{\RdxR} \varphi(\xi) \int_{\R} \theta_{n}^{(t)}(s)\theta(t)dt dT(\xi,s)  = \int_{\RdxR} \varphi(\xi) \theta(s) dT(\xi,s) = \langle T , \varphi \boxtimes \theta \rangle,
\end{equation}
from which we finally obtain
\begin{equation}
\int_{\R} \langle T_{t} , \varphi \rangle \theta(t) dt = \langle T , \varphi \boxtimes \theta \rangle. \quad \blacksquare
\end{equation}

If $T$ is in any of the spaces considered in Theorem \ref{Theo:CadlagRepresentation} and if it satisfies the hypotheses there given, formula \eqref{Eq:Tt-Ts} allows us to obtain an expression for $T_{t} - T_{s}$ for every $t > s$. However, we do not have in general an explicit expression for $T_{t}$ alone, although of course it can be computed following \eqref{Eq:TtAsLimit}. We remark, however, that in the particular cases where $T$ is in $\VRdCfdRpDual$ or in $\MsgRdxRp$, $T_{t} = 0 $ for $t < 0$, and hence it follows immediately that the càdlàg-in-time representation in those cases is given by
\begin{equation}
\langle T_{t} , \varphi \rangle = \langle \frac{\partial T}{\partial t} , \varphi \boxtimes \mathbf{1}_{\left[ 0 , t \right]} \rangle, \quad t \in \Rp,
\end{equation}
for every $\varphi$ in $\VRd$ or in $\CfdRd$ correspondingly.

\section{Analysis of the PDE}
\label{Sec:AnalysisPDE}

Let us now consider the PDE \eqref{Eq:ThePDE}. Let us require $X \in \VRdCfdRDual$. We will see that this requirement will allow us to properly speak about an initial condition for Cauchy problems associated to this equation. In addition, it allows to consider the cases where the spatial symbol function $g$ is continuous, allowing a big variety of pseudo-differential operators. We work in the analogous of a \textit{parabolic framework}, requiring always that $g_{R} \geq 0 $.

\subsection{Solutions over $\RdxRp$}
\label{Sec:SolutionsOverRdxRp}

We begin by looking at for solutions over $\RdxRp$, that is, we do not care if the resulting distributions satisfy the equation over $\RdxRneg$.

Let $X \in \VRdCfdRpDual$. Let us apply the spatial Fourier transform to the equation \eqref{Eq:ThePDE}, obtaining the spatially-multiplicative equation:
\begin{equation}
\label{Eq:TransformedFOEM}
\frac{\partial V}{\partial t} + gV = Y,
\end{equation}
where $Y = \mathscr{F}_{S}(X) \in \MsgRdxRp$ and $V = \mathscr{F}_{S}(U)$ is the transformed unknown. 

We introduce the following operator, which we will call  \textit{Duhamel's Operator}. Consider $\mathcal{D}_{g} : \CfdRdRp \to \CfdRdRp$ defined through
\begin{equation}
\label{Eq:DefDuhamelOperatorTestFunction}
\mathcal{D}_{g}(\psi ) (\xi , t) := \int_{t}^{\infty} e^{-(s-t)g(\xi)} \psi(\xi,s)ds  \quad  (\xi,t) \in \RdxRp.
\end{equation}
We remark that this operator is nothing but a \textit{temporal convolution} (that is, a convolution with respect to the temporal component) with the function $(\xi,s) \in \RdxR \mapsto e^{sg(\xi)}\mathbf{1}_{\Rneg}(s)$ (we may, for instance, extend the domain of $\psi$ to $\RdxR$ by making it null over $\RdxRneg$, in order to properly define such a convolution).

\begin{Prop}
\label{Prop:DgLinearContinuousCfdRdxRp}
Suppose $g : \Rd \mapsto \mathbb{C}$ is a continuous symbol function such that $g_{R} \geq 0 $. Then, $\mathcal{D}_{g}$ is a continuous linear operator from $\CfdRdxRp$ to $\CfdRdxRp$.
\end{Prop}

\textbf{Proof: } The linearity of $\mathcal{D}_{g}$ is straightforward. If we consider a sequence $(\xi_{n} , t_{n} )_{n \in \mathbb{N}} \subset \RdxRp$ such that $(\xi_{n} , t_{n} ) \to (\xi , t ) \in \RdxRp$ as $n \to \infty$, it is immediate to verify using the continuity of $g$ and $\psi$ that
\begin{equation}
\label{Eq:ContinuityAfterDuhamelConvergence}
e^{-(s-t_{n})g(\xi_{n})}\psi(\xi_{n}, s)\mathbf{1}_{\left( t_{n} , \infty \right)}(s) \to e^{-(s-t)g(\xi)}\psi(\xi, s)\mathbf{1}_{\left( t , \infty \right)}(s), \quad \hbox{ as } n \to \infty, \forall s \in \Rp\setminus\lbrace t \rbrace.
\end{equation}
Considering that $g_{R}\geq 0 $, we also have that  $|e^{-ag(\xi)} | \leq 1$ for every $a \geq 0 $, hence
\begin{equation}
\left| e^{-(s-t_{n})g(\xi_{n})}\psi(\xi_{n}, s)\mathbf{1}_{\left( t_{n} , \infty \right)}(s) \right| \leq    \left| \psi(\xi_{n} , s) \right| \leq \frac{|\psi(\xi_{n},s)|}{(1+s^{2})}(1+s^{2}) \leq \frac{\sup_{(\eta,u) \in \RdxRp}\left| \psi( \eta , u)(1+u^{2})  \right|}{1+s^{2}}. 
\end{equation}
We conclude that the convergence \eqref{Eq:ContinuityAfterDuhamelConvergence} is dominated by  function $ s \mapsto \frac{\sup_{(\eta,u) \in \RdxRp}\left| \psi( \eta , s)(1+s^{2})  \right|}{1+s^{2}}, $ which is integrable over $\Rp$ since $\psi \in C_{FD}(\RdxRp)$. It follows from Dominated Convergence Theorem that the function $\mathcal{D}_{g}(\psi)$ is continuous. Let $n_{S},n_{T} \in \mathbb{N}$. We consider that
\small
\begin{equation}
\label{Eq:DgLinearConitnuousRdxRp}
\left.\begin{aligned}
\sup_{(\xi,t) \in \RdxRp}\left| (1+|\xi|^{2})^{n_{S}}(1+t^{2})^{n_{T}}\mathcal{D}_{g}(\psi) (\xi,t) \right| &= \sup_{ (\xi,t) \in \RdxRp } \left|  (1+|\xi|^{2})^{n_{S}}(1+t^{2})^{n_{T}}\int_{t}^{\infty} e^{-(s-t)g(\xi)}\psi(\xi,s)ds   \right|  \\
&\leq \sup_{(\xi,t) \in \RdxRp} \Big\lbrace \int_{t}^{\infty}(1+|\xi|^{2})^{n_{S}}(1+t^{2})^{n_{T}} |\psi(\xi,s)|ds  \Big\rbrace  \\
&\leq \sup_{(\xi,t) \in \RdxRp} \Big\lbrace \int_{t}^{\infty}(1+|\xi|^{2})^{n_{S}}(1+s^{2})^{n_{T}} |\psi(\xi,s)|ds  \Big\rbrace  \\
&\leq \sup_{(\xi,t) \in \RdxRp} \Big\lbrace \int_{\Rp}(1+|\xi|^{2})^{n_{S}}(1+s^{2})^{n_{T}+1} |\psi(\xi,s)|\frac{ds}{(1+s^{2})}  \Big\rbrace  \\
 &\leq  \frac{\pi}{4} \sup_{(\xi,s) \in \RdxRp} \left|(1+|\xi|^{2})^{n_{S}}(1+s^{2})^{n_{T}+1} \psi(\xi,s)\right|.
\end{aligned}\right.
\end{equation}
\normalsize
This proves that $\mathcal{D}_{g}(\psi)$ is a fast-decreasing function. In addition, this also proves that $\mathcal{D}_{g}$ is a continuous linear operator from $C_{FD}(\RdxRp)$ to $C_{FD}(\RdxRp)$. $\blacksquare $

The adjoint operator of Duhamel's operator is denoted by $\mathcal{D}_{g}^{*}$ and it is hence a continuous linear operator from $\mathscr{M}_{SG}(\RdxRp)$ to $\mathscr{M}_{SG}(\RdxRp)$.

\begin{Prop}
\label{Prop:D*gYsatisfiesPDE}
Let $Y \in  \mathscr{M}_{SG}(\RdxRp)\subset \TemperedRdxR$. Then, $\mathcal{D}_{g}^{*}(Y)$ satisfies \eqref{Eq:TransformedFOEM} in the sense of $\TemperedRdxR$.
\end{Prop}

\textbf{Proof: } Let $\psi \in \SchwartzRdxR$. We have
\begin{equation}
\langle \frac{\partial \mathcal{D}_{g}^{*}(Y)}{\partial t} , \psi \rangle = - \langle \mathcal{D}_{g}^{*}(Y) , \frac{\partial \psi}{\partial t} \rangle = - \langle \mathcal{D}_{g}^{*}(Y) , \frac{\partial \psi}{\partial t}\big|_{\RdxRp} \rangle = - \langle Y , \mathcal{D}_{g}\left(\frac{\partial \psi}{\partial t}\big|_{\RdxRp}\right) \rangle.
\end{equation} 
By integrations by parts, for every $(\xi,t) \in \RdxRp$ it holds
\begin{equation}
\label{Eq:-DgDerPsi}
\left. \begin{aligned}
- \mathcal{D}_{g}\left(\frac{\partial \psi}{\partial t}\Big|_{\RdxRp}   \right)(\xi,t) &= -\int_{t}^{\infty} e^{-(s-t)g(\xi)} \frac{\partial \psi}{\partial t}(\xi,s) ds \\
&= - \left[  e^{-(s-t)g(\xi)}\psi(\xi,s) \big|^{s=\infty}_{s=t} - \int_{t}^{\infty} e^{-(s-t)g(\xi)}(-g(\xi))\psi(\xi,s)ds    \right] \\
&= - \left[ - \psi(\xi,t) + g(\xi)\int_{t}^{\infty} e^{-(s-t)g(\xi)}\psi(\xi,s)ds \right] \\
&= \psi(\xi,t) - g(\xi) \mathcal{D}_{g}(\psi)(\xi,t).
\end{aligned} \right.
\end{equation}
It is immediate that $g\mathcal{D}_{g}(\psi) = \mathcal{D}_{g}( g\psi ) $. We obtain thus
\begin{equation}
\langle \frac{\partial \mathcal{D}_{g}^{*}(Y)}{\partial t} , \psi \rangle =  \langle Y , \psi \rangle - \langle Y , \mathcal{D}_{g}( g\psi ) \rangle = \langle Y , \psi \rangle - \langle \mathcal{D}_{g}^{*}(Y) , g\psi \rangle = \langle Y , \psi \rangle - \langle g\mathcal{D}_{g}^{*}(Y) , \psi \rangle,
\end{equation}
where the equality $\langle \mathcal{D}_{g}^{*}(Y) , g\psi \rangle = \langle g\mathcal{D}_{g}^{*}(Y) , \psi \rangle $ is justified since $\mathcal{D}_{g}^{*}(Y)$ is a slow-growing measure and $g$ is a polynomially bounded continuous function. This proves that $\mathcal{D}_{g}^{*}(Y)$ satisfies the PDE \eqref{Eq:TransformedFOEM} in the sense of $\TemperedRdxR$. $\blacksquare$

\begin{Corol}
\label{Corol:DerivativeDuhamelInMsgRdRp}
Let $Y \in \mathscr{M}_{SG}(\RdxRp)$. Then, $ \frac{\partial \mathcal{D}_{g}^{*}(Y)}{\partial t} \in \mathscr{M}_{SG}(\RdxRp) $. 
\end{Corol}

\textbf{Proof: } $\frac{\partial \mathcal{D}_{g}^{*}(Y)}{\partial t} = Y - g\mathcal{D}_{g}^{*}(Y) \in \mathscr{M}_{SG}(\RdxRp) $. $\blacksquare$

From Theorem \ref{Theo:CadlagRepresentation}, the next Corollary follows immediately.

\begin{Corol}
\label{Corol:DuhammelTHasCadlagRepresentation}
Let $Y \in \mathscr{M}_{SG}(\RdxRp)$. Then, $\mathcal{D}_{g}^{*}(Y)$ has a càdlàg-in-time representation.
\end{Corol}

Let us describe the càdlàg-in-time representation of $\mathcal{D}_{g}^{*}(Y)$ for any $Y \in \MsgRdxRp$, which we denote $(\mathcal{D}_{g}^{*}(Y)_{t} )_{t \in \Rp}$. We consider $t \in \Rp$, $(\theta_{n}^{(t)})_{n \in \mathbb{N}} \subset \DR$ a sequence such that $ \theta_{n}^{(t)} \to \delta_{t}^{+}$, for which we suppose $\supp(\theta_{n}^{(t)}) \subset \left[ t , t+1 \right] $ for every $n \in \mathbb{N}$. Let $\varphi \in \CfdRd$. We have thus,
\begin{equation}
\left. \begin{aligned}
\langle \mathcal{D}_{g}^{*}(Y) , \varphi \boxtimes \theta_{n}^{(t)} \rangle &= \langle Y , \mathcal{D}_{g}(  \varphi \boxtimes \theta_{n}^{(t)}   ) \rangle \\
&= \int_{\RdxRp} \mathcal{D}_{g}(\varphi\boxtimes \theta_{n}^{(t)})(\xi,s) dY(\xi,s) \\
&= \int_{\RdxRp} \int_{\Rp} e^{-(u-s)g(\xi)}\mathbf{1}_{\left[ s , \infty \right)}(u) \varphi(\xi)\theta_{n}^{(t)}(u)du dY(\xi,s).
\end{aligned} \right.
\end{equation}
Since $\theta_{n}^{(t)} \to \delta_{t}^{+}$, and since the function $u \mapsto e^{-(u-s)g(\xi)}\mathbf{1}_{\left[ s , \infty \right)}(u)\varphi(\xi)$ is càdlàg for every $(\xi,s) \in \RdxRp$, one has
\begin{equation}
\int_{\Rp} e^{-(u-s)g(\xi)}\mathbf{1}_{\left[ s , \infty \right)}(u) \varphi(\xi)\theta_{n}^{(t)}(u)du \to \varphi(\xi) e^{-(t-s)g(\xi)}\mathbf{1}_{\left[ s , \infty \right)}(t), \quad \hbox{ as }  n \to \infty, \forall (\xi,s) \in \RdxRp.
\end{equation}
In addition, by playing with the supports of the functions involved, one obtains
\begin{equation}
\left|  \int_{\Rp} e^{-(u-s)g(\xi)}\mathbf{1}_{\left[ s , \infty \right)}(u) \varphi(\xi)\theta_{n}^{(t)}(u)du \right| \leq |\varphi(\xi)| \mathbf{1}_{\left[0 , t+1 \right]}(s), \quad \forall (\xi,s) \in \RdxRp, \forall n \in \mathbb{N}.
\end{equation}
Since $(\xi,s) \mapsto   |\varphi(\xi)| \mathbf{1}_{\left[0 , t+1 \right]}(s)$ is integrable with respect to $Y$ over $\RdxRp$, one obtains from Dominated Convergence Theorem
\begin{equation}
\lim_{n \to \infty} \langle \mathcal{D}_{g}^{*}(Y) , \varphi\boxtimes \theta_{n}^{(t)} \rangle = \int_{\Rd\times\Rp} \varphi(\xi) e^{-(t-s)g(\xi)}\mathbf{1}_{\left[ s , \infty \right)}(t) dY(\xi,s) = \int_{\Rd\times\left[0 , t \right]} \varphi(\xi) e^{-(t-s)g(\xi)} dY(\xi,s),
\end{equation}
from which we conclude
\begin{equation}
\label{Eq:CadlagRepresentationDuhamelRdxRp}
\langle \mathcal{D}_{g}^{*}(Y)_{t} , \varphi \rangle = \int_{\Rd\times\left[0,t\right]} \varphi(\xi) e^{-(t-s)g(\xi)}dY(\xi,s), \quad \forall \varphi \in \CfdRd, \forall t \in \Rp.
\end{equation}

We conclude that we can always find a solution to the transformed problem \eqref{Eq:TransformedFOEM} which has a càdlàg-in-time behaviour, and hence for which the notion of an \textit{initial condition} makes sense. Nevertheless, the following result shows that we still have some difficulties if we want to consider any arbitrary initial condition in a Cauchy problem associated to equation \eqref{Eq:TransformedFOEM}.

\begin{Prop}
\label{Prop:Dg*YUniqueSolutionMsgRdRp}
$\mathcal{D}_{g}^{*}(Y)$ is the unique possible solution in $\mathscr{M}_{SG}(\RdxRp)$ to equation \eqref{Eq:TransformedFOEM}. 
\end{Prop}

This result follows from a simple fact, which is actually an equivalent statement to Proposition \ref{Prop:Dg*YUniqueSolutionMsgRdRp}: \textit{the homogeneous problem}
\begin{equation}
\label{Eq:HomogenoeusTransformedFOEM}
\frac{\partial V_{H}}{\partial t} + gV_{H} = 0
\end{equation}
\textit{has no non-trivial solutions in } $\mathscr{M}_{SG}(\RdxRp)$.

\textbf{Proof:}  Let us suppose there are two solutions in $\mathscr{M}_{SG}(\RdxRp)$, say $V_{1}$ and $V_{2}$. Then, by linearity of the equation, the difference $V_{H} = V_{1}-V_{2} \in \mathscr{M}_{SG}(\RdxRp)$ must satisfy the homogeneous problem \eqref{Eq:HomogenoeusTransformedFOEM}. Let us look for solutions to the homogeneous problem in the bigger space of measures $\mathscr{M}(\RdxR)$. With a typical analysis we have that $V_{H}$ satisfies
\begin{equation}
\frac{\partial}{\partial t} \left( e^{tg}V_{H} \right) = 0.
\end{equation}
Hence,
\begin{equation}
\label{Eq:etgVh=Sboxtimes1}
e^{tg}V_{H} = S \boxtimes \mathbf{1}
\end{equation}
for some $S \in \mathscr{D}'(\Rd)$, and since we have required that $V_{H}$ must be a measure, $S$ must be in $\mathscr{M}(\Rd)$. It turns out that $V_{H}$ is of the form
\begin{equation}
\label{Eq:ExpressionVHomogenoeusTransformedFOEM}
V_{H} = e^{-tg}\left( S \boxtimes \mathbf{1} \right).
\end{equation}
However, expression \eqref{Eq:ExpressionVHomogenoeusTransformedFOEM} does not provide a measure with support on $\RdxRp$ unless $S = 0$. Hence, there is no solution to \eqref{Eq:HomogenoeusTransformedFOEM} in $\mathscr{M}_{SG}(\RdxRp)$ besides the trivial solution.  $\blacksquare$

In the last proof we remark that if we consider the \textit{restriction} of the measure $V_{H}$ over $\RdxRp$,  $\mathbf{1}_{\RdxRp}V_{H}$, then we do obtain a measure in $\mathscr{M}_{SG}(\RdxRp)$ if $S \in \mathscr{M}_{SG}(\Rd)$. However, this measure does not satisfy \eqref{Eq:HomogenoeusTransformedFOEM} in the sense of $\TemperedRdxR$. It rather satisfies
\begin{equation}
\frac{\partial  }{\partial t}(\mathbf{1}_{\RdxRp}V_{H}) + g\mathbf{1}_{\RdxRp}V_{H} = S \boxtimes \delta_{0_{T}},
\end{equation}
where $\delta_{0_{T}} \in \MsgRp$ denotes the Dirac measure at $0 \in \Rp$.

\begin{Corol}
There exists a unique solution in $\VRdCfdRpDual$ to equation \eqref{Eq:ThePDE}.
\end{Corol}
\textbf{Proof: } Take $U = \mathscr{F}_{S}^{-1}\left(\mathcal{D}_{g}^{*}(Y)\right) = \mathscr{F}_{S}^{-1}\left(\mathcal{D}_{g}^{*}(\mathscr{F}_{S}(X))\right)$. $\blacksquare$

\subsection{Solutions over $\RdxR$}
\label{Sec:SolutionsOverRdxR}

Consider for now the PDE \eqref{Eq:ThePDE} with the only condition $X \in \VRdxRDual$. Let us make an extra supposition for $g$: there exists $\kappa > 0 $ such that $g_{R} \geq \kappa $. By following Section \ref{Sec:OperatorsSymbol}, we can guarantee the existence of a unique solution in $\VRdxRDual$, which we denote by $U^{\infty}$ and which is given by
\begin{equation}
\label{Eq:Uinfinity}
U^{\infty} =  \mathscr{F}^{-1}\left(  \frac{1}{i\omega + g(\xi)}\mathscr{F}(X) \right).
\end{equation}
In general, it is not clear if $U^{\infty}$ has a càdlàg-in-time representation. We cannot simply apply the principles of Theorem \ref{Theo:CadlagRepresentation} since the distributions in $\VRdxRDual$ do not necessarily behave as a measure in time. In order to restrain our work to such kind of solutions, let us consider the case where $X \in \VRdxRDual \cap  \VRdCfdRDual$. In such case $X$ is a measure in the temporal component and both $\mathscr{F}_{S}(X)$ and $\mathscr{F}(X)$ are slow-growing measures over $\RdxR$.

Consider the spatially transformed problem \eqref{Eq:TransformedFOEM}. Let us reconsider Duhamel's operator, defined as in \eqref{Eq:DefDuhamelOperatorTestFunction}, but this time we suppose $\psi \in \CfdRdxR$, with $(\xi,t) \in \RdxR$.

\begin{Prop}
\label{Prop:DgLinearContinuousCfdRdxR}
Suppose $g : \Rd \mapsto \mathbb{C}$ is a continuous symbol function such that $g_{R} \geq \kappa$ for some $ \kappa > 0 $. Then, $\mathcal{D}_{g}$ defines a continuous linear operator from $\CfdRdxR$ to $\CfdRdxR$.
\end{Prop}
We will use the following simple Lemma:
\begin{Lemma}
\label{Lemma:TraceCfdRdxRmInCfdRd}
Let $\psi \in C_{FD}(\RdxRm)$, with $d,m \in \mathbb{N}_{*}$. Then, $x \in \Rd \mapsto \sup_{y \in \Rm}\left| (1+|y|^{2})^{M}\psi(x,y) \right|$ is in $\CfdRd$ for every $M \in \mathbb{N}$.
\end{Lemma}
\textbf{Proof: } The proof of the continuity follows immediately from the uniform continuity of the functions in $C_{FD}(\RdxRm)$. The proof of the fast-decreasing behaviour follows from $$ (1+|x|^{2})^{N}\sup_{y \in \Rm}\left| (1+|y|^{2})^{M}\psi(x,y) \right|  \leq \sup_{(x,y) \in \RdxRm} \left|(1+|x|^{2})^{N} (1+|y|^{2})^{M}\psi(x,y) \right| < \infty,  $$
for all $x \in \Rd$ and for all $N \in \mathbb{N}$. $\blacksquare$

\textbf{Proof of Proposition \ref{Prop:DgLinearContinuousCfdRdxR}:} The linearity of $\mathcal{D}_{g}$ is straightforward. Let $\psi \in \CfdRdxR$. The continuity of $\mathcal{D}_{g}(\psi)$ can be proven following the same arguments as in Proposition \ref{Prop:DgLinearContinuousCfdRdxRp}. Let $n_{T}, n_{S} \in \mathbb{N}$. From Lemma \ref{Lemma:ConvexityInequalities}, we obtain for any $(\xi,t) \in \RdxR$, 
\small
\begin{equation}
\left. \begin{aligned}
\left| (1+|\xi|^{2})^{n_{S}}(1+t^{2})^{n_{T}} \mathcal{D}_{g}(\psi) \right| &= \left|  (1+|\xi|^{2})^{n_{S}}(1+t^{2})^{n_{T}} \int_{t}^{\infty} e^{-(s-t)g(\xi)} \psi(\xi,s)ds  \right| \\
&\leq \int_{t}^{\infty}e^{-(s-t)\kappa} (1+t^{2})^{n_{T}}(1+|\xi|^{2})^{n_{S}} |\psi(\xi,s)| ds \\
&\leq 2^{n_{T}-1}\int_{t}^{\infty}e^{-(s-t)\kappa}\left( (1+2(s-t)^{2})^{n_{T}} + 2^{n_{T}}s^{2n_{T}} \right) \sup_{\xi \in \Rd} \left| (1+|\xi|^{2})^{n_{S}}\psi(\xi,s) \right| ds \\
&\leq 2^{n_{T}-1} \int_{t}^{\infty}e^{-(s-t)\kappa}(1+2(s-t)^{2})^{n_{T}}\sup_{\xi \in \Rd} \left| (1+|\xi|^{2})^{n_{S}}\psi(\xi,s) \right| ds  \\
&\quad + 2^{2n_{T} - 1}\int_{t}^{\infty}e^{-(s-t)\kappa}s^{2n_{T}} \sup_{\xi \in \Rd} \left| (1+|\xi|^{2})^{n_{S}}\psi(\xi,s) \right|ds. 
\end{aligned}\right. 
\end{equation}
\normalsize
From Lemma \ref{Lemma:TraceCfdRdxRmInCfdRd}, the function $s \in \R \mapsto \sup_{\xi \in \Rd} \left| (1+|\xi|^{2})^{n_{S}}\psi(\xi,s) \right| $ is in $\CfdR$ hence it is bounded. In addition, the function $ s \in \R \mapsto e^{s\kappa}(1+2s^{2})^{n_{T}}\mathbf{1}_{\Rneg}(s) $ is integrable over $\R$. Using that   $\| f_{1} \ast f_{2} \|_{L^{\infty}} \leq \| f_{2} \|_{L^{1}} \| f_{1} \|_{L^{\infty}}$ for $f_{1} \in L^{\infty}(\R)$ and $f_{2} \in L^{1}(\R)$, one has
\small
\begin{equation}
\int_{t}^{\infty}e^{-(s-t)\kappa}(1+2(s-t)^{2})^{n_{T}}\sup_{\xi \in \Rd} \left| (1+|\xi|^{2})^{n_{S}}\psi(\xi,s) \right| ds \leq \int_{\Rneg} e^{s\kappa}(1+2s^{2})^{n_{T}} ds \sup_{(\xi,t) \in \RdxR}\left| (1+|\xi|^{2})^{n_{S}}\psi(\xi , t) \right|. 
\end{equation} 
\normalsize
On the other hand, one has
\begin{equation}
\left. \begin{aligned}
\int_{t}^{\infty}e^{-(s-t)\kappa}s^{2n_{T}} \sup_{\xi \in \Rd} \left| (1+|\xi|^{2})^{n_{S}}\psi(\xi,s) \right|ds &\leq \int_{t}^{\infty}(1+s^{2})^{n_{T}+1} \sup_{\xi \in \Rd} \left| (1+|\xi|^{2})^{n_{S}}\psi(\xi,s) \right|\frac{ds}{1+s^{2}}   \\
&\leq \frac{\pi}{2} \sup_{(\xi,t) \in \RdxR}\left| (1+|\xi|^{2})^{n_{S}}(1+t^{2})^{n_{T}+1}\psi(\xi , t) \right|. 
\end{aligned} \right.
\end{equation}
We conclude thus
\begin{equation}
\label{Eq:DgLinearContinuousCfdRdxR}
\left. \begin{aligned}
\sup_{(\xi,t) \in \RdxR}\Bigg| &(1+|\xi|^{2})^{n_{S}}(1+t^{2})^{n_{T}} \mathcal{D}_{g}(\psi)(\xi,t) \Bigg|  \\
&\leq 2^{n_{T}-1}\left( \int_{\Rneg}e^{s\kappa}(1+2s^{2})^{n_{T}}ds + \pi 2^{n_{T}-2} \right)\sup_{(\xi,t) \in \RdxR }\left| (1+|\xi|^{2})^{n_{S}}(1+t^{2})^{n_{T}+1}\psi(\xi,t)  \right|.
\end{aligned} \right.
\end{equation}
This proves both that $\mathcal{D}_{g}(\psi)$ is in $\CfdRdxR$ and that $\mathcal{D}_{g}$ defines a continuous linear operator from $\CfdRdxR$ to $\CfdRdxR$. $\blacksquare$

The adjoint of Duhamel's operator, $\mathcal{D}_{g}^{*}$, is then a continuous linear operator from $\MsgRdxR$ to $\MsgRdxR$. In addition, we have the following simple manner of describing such operator.

\begin{Prop}
\label{Prop:DgEqualsFTdivisionFTinverse}
$\mathcal{D}_{g} = \mathscr{F}_{T}\left( \frac{1}{i\omega + g}\mathscr{F}_{T}^{-1}\left( \cdot \right) \right)$ over $\CfdRdxR$.
\end{Prop}

\textbf{Proof: } Let $\psi \in \CfdRdxR$. Since $g_{R} \geq \kappa > 0 $, the function $\omega \in \R \mapsto \frac{1}{i\omega + g(\xi)}$ is in $\OMR$ for every $\xi \in \Rd$. Hence, applying the exchange formula for the temporal Fourier transform \eqref{Eq:ExchangeFormula}, we have 
\begin{equation}
\mathscr{F}_{T}\left( (\xi,\omega) \mapsto \frac{1}{i\omega + g(\xi)} \mathscr{F}_{T}^{-1}(\psi)(\xi, \omega) \right) = \frac{1}{\sqrt{2\pi}} \mathscr{F}_{T}\left( \frac{1}{i\omega + g} \right) \stackrel{(\R)}{\ast} \psi,
\end{equation}
where $\stackrel{(\R)}{\ast}$ denotes a temporal convolution. A typical result from Fourier Analysis, which can be applied considering that $g_{R} \geq \kappa > 0 $, states that
\begin{equation}
\mathscr{F}_{T}\left(  \frac{1}{i\omega + g} \right)(\xi,t) = \sqrt{2\pi}e^{tg(\xi)}\mathbf{1}_{\mathbb{R}_{*}^{-}}(t)
\end{equation}
in distributional sense. We conclude that $\mathscr{F}_{T}\left( \frac{1}{i\omega + g}\mathscr{F}_{T}^{-1}(\cdot) \right) $ equals the temporal convolution with the function $(\xi,s) \mapsto e^{sg(\xi)}\mathbf{1}_{\Rneg}(s)$, which is no other than Duhamel's operator over $\CfdRdxR$. $\blacksquare$

\begin{Corol}
\label{Corol:Dg*EqualsFTinverseDivisionFT}
$\mathcal{D}_{g}^{*}(\mu) = \mathscr{F}_{T}^{-1}\left( \frac{1}{i\omega + g}\mathscr{F}_{T}(\mu) \right)$ for every $\mu \in \MsgRdxR$.
\end{Corol}

\textbf{Proof: } Since the adjoint of the operator $\mathscr{F}_{T}\left( \frac{1}{i\omega + g}\mathscr{F}_{T}^{-1}\left( \cdot \right) \right)$ is $\mathscr{F}_{T}^{-1}\left( \frac{1}{i\omega + g}\mathscr{F}_{T}\left( \cdot \right) \right)$, the result follows from Riesz Representation Theorem \ref{Theo:RieszReperesentationDualCfdMsg}. $\blacksquare$

We can now conclude the following result: when $X \in \VRdxRDual\cap\VRdCfdRDual$, the solution $U^{\infty}$ given by \eqref{Eq:Uinfinity} can be expressed through
\begin{equation}
U^{\infty} = \mathscr{F}^{-1}\left( \frac{1}{i\omega + g} \mathscr{F}(X) \right) = \mathscr{F}_{S}^{-1}\left( \mathscr{F}_{T}^{-1}\left( \frac{1}{i\omega + g} \mathscr{F}_{T}\left(\mathscr{F}_{S}(X) \right) \right) \right) = \mathscr{F}_{S}^{-1}\left( \mathcal{D}_{g}^{*}( \mathscr{F}_{S}(X) ) \right). 
\end{equation}  
We remark that this implies immediately that $U^{\infty} \in \VRdxRDual \cap \VRdCfdRDual$.

\begin{Prop}
\label{Prop:UinfinityHasCadlagRepresentation}
$U^{\infty}$ has a càdlàg-in-time representation.
\end{Prop}
\textbf{Proof: } Let $V^{\infty} = \mathscr{F}_{S}(U^{\infty}) = \mathcal{D}_{g}^{*}(\mathscr{F}_{S}(X)) \in \MsgRdxR$. $V^{\infty}$ satisfies
\begin{equation}
\frac{\partial V^{\infty}}{\partial t} = -gV^{\infty} + \mathscr{F}_{S}(X) = -g \mathcal{D}_{g}^{*}(\mathscr{F}_{S}(X)) - \mathscr{F}_{S}(X).
\end{equation}
Since $-g \mathcal{D}_{g}^{*}(\mathscr{F}_{S}(X)) - \mathscr{F}_{S}(X) \in \MsgRdxR$, from Theorem \ref{Theo:CadlagRepresentation}, $V^{\infty}$ has a càdlàg-in-time representation and so does $U^{\infty}$ from Proposition \ref{Prop:CadlagRepresentationInterchangeSpatialFourier}. $\blacksquare$

An explicit expression for the càdlàg-in-time representation of $\mathscr{F}_{S}(U^{\infty})$ can be obtained applying the same procedure used to obtain \eqref{Eq:CadlagRepresentationDuhamelRdxRp} in Section \ref{Sec:SolutionsOverRdxRp}. The càdlàg-in-time representation of $U^{\infty}$, denoted $(U^{\infty}_{t})_{t \in \R} \subset \VRdDual$, is obtained by applying and inverse spatial Fourier transform. The details are left to the reader. The final expression for the càdlàg-in-time representation of $U^{\infty}$ is given by
\begin{equation}
\label{Eq:CadlagRepresentationUinfinity}
\langle U^{\infty}_{t} , \varphi \rangle = \int_{\Rd\times\left( -\infty , t \right]} e^{-(t-s)g(\xi)} \mathscr{F}_{S}^{-1}(\varphi)(\xi) d\mathscr{F}_{S}(X)(\xi,s), \quad \forall \varphi \in \VRd, \forall t \in \R.
\end{equation}

\section{The Cauchy problem}
\label{Sec:CauchyProblem}

Let us right now focus on the Cauchy problem
\begin{equation}
\label{Eq:CauchyProblemFOEMFormal}
\left\lbrace\begin{array}{lcr}
\dfrac{\partial U}{\partial t} + \mathcal{L}_{g}U = X & & \hbox{over} \  \RdxRp \\
U \big|_{t=0} = U_{0} & &
\end{array},
\right.
\end{equation}
for suitable distributions $X$ and $U_{0}$. Let us require $ X \in \VRdCfdRpDual$ and $U_{0} \in \VRdDual$. We apply the spatial Fourier transform and we work with the transformed Cauchy problem
\begin{equation}
\label{Eq:CauchyProblemTransformedFOEMFormal}
\left\lbrace\begin{array}{lcc}
\dfrac{\partial V}{\partial t} + gV = Y & & \\
V \big|_{t=0} = V_{0} & &
\end{array},
\right.
\end{equation}
where $Y = \mathscr{F}_{S}(X) \in \mathscr{M}_{SG}(\RdxRp)$, $V_{0} = \mathscr{F}_{S}(U_{0}) \in \mathscr{M}_{SG}(\Rd)$ and $V = \mathscr{F}_{S}(U)$ is the transformed unknown. 

Proposition \ref{Prop:Dg*YUniqueSolutionMsgRdRp} guarantees that if we look at for solutions in $\mathscr{M}_{SG}(\RdxRp)$ to the Cauchy problem \eqref{Eq:CauchyProblemTransformedFOEMFormal}, we will have just one possibility: $\mathcal{D}_{g}^{*}(Y)$, which has a càdlàg-in-time representation. $\mathcal{D}_{g}^{*}(Y)$ is a solution to the Cauchy problem \eqref{Eq:CauchyProblemTransformedFOEMFormal}  if and only if  $V_{0} = \mathcal{D}_{g}^{*}(Y)_{0} $. Hence, we have no freedom at all to fix an arbitrary initial condition. The next result proves that, with some extra arrangements, we can gain more freedom in the initial condition by requiring that the restriction to $\RdxRp$ of the solution is in $\mathscr{M}_{SG}(\RdxRp)$, rather than the solution itself.

\begin{Theo}
\label{Theo:ExistenceUniquenessSolutionTransformedFOEM}
Let $Y \in \mathscr{M}_{SG}(\RdxRp)$. Let $g : \Rd \to \mathbb{C}$ be a continuous spatial symbol function such that $g_{R} \geq 0 $. Let $V_{0} \in \mathscr{M}_{SG}(\Rd)$. Then, there exists a unique measure $V \in \mathscr{M}(\RdxR)$ such that
\begin{itemize}
\item $\frac{\partial V}{\partial t} + gV = Y$ in the sense of $\mathscr{D}'(\RdxR)$,
\item its restriction to $\RdxRp$ is in $\mathscr{M}_{SG}(\RdxRp)$ and it has a càdlàg-in-time representation whose evaluation at $t=0$ is $V_{0}$.
\end{itemize}
\end{Theo}

\textbf{Proof: } For the existence, consider the measure $V \in \mathscr{M}(\RdxR)$ defined by
\begin{equation}
\label{Eq:ExpressionSolutionVTransformedFOEMugly}
V = e^{-tg}\left( (V_{0} - \mathcal{D}_{g}^{*}(Y)_{0})\boxtimes \mathbf{1} \right) + \mathcal{D}_{g}^{*}(Y),
\end{equation}
which is more explicitly expressed as (we recall that $\mathcal{D}_{g}^{*}(Y)$ is a measure with support in $\RdxRp$ and $V_{0}$ and $\mathcal{D}_{g}^{*}(Y)_{0}$ are measures over $\Rd$):
\begin{equation}
\label{Eq:ExplicitVTestFunctionUgly}
\langle V , \psi \rangle = \int_{\mathbb{R}}\int_{\Rd} e^{-tg(\xi)}\psi(\xi,t) d(   V_{0} - \mathcal{D}_{g}^{*}(Y)_{0})(\xi)dt + \int_{\RdxRp} \psi(\xi,t) d \mathcal{D}_{g}^{*}(Y) (\xi,t), \quad \forall \psi \in \mathscr{D}(\RdxR).
\end{equation}
The fact that $V$ satisfies the PDE \eqref{Eq:TransformedFOEM} arises from Proposition \ref{Prop:D*gYsatisfiesPDE} and from the fact that the measure $ e^{-tg}\left( (V_{0} - \mathcal{D}_{g}^{*}(Y)_{0})\boxtimes \mathbf{1} \right)$ satisfies the homogeneous equation \eqref{Eq:HomogenoeusTransformedFOEM}, since it is of the form \eqref{Eq:etgVh=Sboxtimes1}.

The restriction of $V$ over $\RdxRp$ is given by 
\begin{equation}
\label{Eq:RestrictionSolutionTransformedCauchyFOEM}
\mathbf{1}_{\RdxRp}V = e^{-tg}\left( (V_{0} - \mathcal{D}_{g}^{*}(Y)_{0})\boxtimes \mathbf{1}_{\Rp} \right) + \mathcal{D}_{g}^{*}(Y).
\end{equation}
It is immediate that $\mathbf{1}_{\RdxRp}V$ is in $\mathscr{M}_{SG}(\RdxRp)$ since $g_{R} \geq 0 $. A typical computation using the derivative of the product and Proposition \ref{Prop:D*gYsatisfiesPDE} allows to conclude that
\begin{equation}
\label{Eq:TemporalDerivvative1RdxRpV}
\frac{\partial}{\partial t}\left( \mathbf{1}_{\RdxRp}V \right) = -ge^{-tg}\left( V_{0} - \mathcal{D}_{g}^{*}(Y)  \right)\boxtimes \mathbf{1}_{\Rp} + e^{-gt}\left( V_{0} - \mathcal{D}_{g}^{*}(Y)_{0} \right)\boxtimes \delta_{0_{T}} - g\mathcal{D}_{g}^{*}(Y) + Y, 
\end{equation}
where $\delta_{0_{T}}$ is the Dirac measure at $0 \in \Rp$. We remark that $$e^{-gt}\left( V_{0} - \mathcal{D}_{g}^{*}(Y)_{0} \right)\boxtimes \delta_{0_{T}} = \left( V_{0} - \mathcal{D}_{g}^{*}(Y)_{0} \right)\boxtimes \delta_{0_{T}}.$$ From this we obtain that the temporal derivative of $\mathbf{1}_{\RdxRp}V$ is a measure in $\mathscr{M}_{SG}(\RdxRp)$ since it is a sum of measures in $\mathscr{M}_{SG}(\RdxRp)$. From Theorem \ref{Theo:CadlagRepresentation} it follows that $\mathbf{1}_{\RdxRp}V$ has a càdlàg-in-time representation. One can apply formula \eqref{Eq:TtAsLimit} to compute the càdlàg-in-time representation of $V$ over $\RdxRp$ (which is evidently the same as the representation of $\mathbf{1}_{\RdxRp}V$), obtaining
\begin{equation}
\label{Eq:CadlagRepresentationOfV}
V_{t} = e^{-tg}\left( V_{0} - \mathcal{D}_{g}^{*}(Y)_{0} \right) + \mathcal{D}_{g}^{*}(Y)_{t} \in \mathscr{M}_{SG}(\Rd), \quad \forall t \in \Rp.
\end{equation}
And from this, it is immediate that the evaluation at $t=0$ if this càdlàg representation equals the desired initial condition $V_{0}$. This proves the existence.

If we suppose that there are two measures $V_{1}$ and $V_{2}$ satisfying the conditions in Theorem \ref{Theo:ExistenceUniquenessSolutionTransformedFOEM}, we consider then the difference $V_{H} = V_{1} - V_{2}$ must satisfies the homogenoeus problem \eqref{Eq:HomogenoeusTransformedFOEM}, and hence it must be of the form \eqref{Eq:etgVh=Sboxtimes1} for some $S \in \mathscr{M}_{SG}(\Rd)$. But this implies that the evaluation at $0$ of its càdlàg-in-time representation is $V_{H,0} = S $. Since in addition, $V_{1}\big|_{t=0} = V_{2}\big|_{t=0}$, then it follows that $V_{H,0}$ must be null, and hence $S=0$. It follows that $V_{H} = 0$. This proves that $V$ is the unique solution to \eqref{Eq:CauchyProblemTransformedFOEMFormal} satisfying the required properties. $\blacksquare$

As stated in the proof of Theorem \ref{Theo:ExistenceUniquenessSolutionTransformedFOEM}, the solution $V$ has a càdlàg-in-time representation over $\Rp$ given by \eqref{Eq:CadlagRepresentationOfV}. Using the càdlàg-in-time representation of Duhamel's operator over $\MsgRdxRp$  \eqref{Eq:CadlagRepresentationDuhamelRdxRp}, one obtains that this family is given by
\begin{equation}
\label{Eq:CadlagRepresentationSolutionTransformedCauchyProblem}
\langle V_{t} , \varphi \rangle = \int_{\Rd} e^{-tg(\xi)} \varphi(\xi) dV_{0}(\xi) + \int_{\Rd \times\left( 0 , t \right]} e^{-(t-s)g(\xi)}\varphi(\xi)dY(\xi,s), \quad \forall \varphi \in \CfdRd, \forall t \in \Rp.
\end{equation}

\begin{Theo}
\label{Theo:ExistenceUniquenessSolutionCuachyFOEMformal}
Let $X \in \VRdCfdRpDual$. Let $g : \Rd \to \mathbb{C}$ be a continuous symbol function such that $g_{R} \geq 0 $. Let $U_{0} \in \VRdDual$. Then, there exists a unique distribution $U \in \VRdCfdRpDual$ such that
\begin{itemize}
\item It has a càdlàg-in-time representation whose evaluation at $t=0$ is $U_{0}$.
\item It satisfies
\begin{equation}
\label{Eq:UsolvesPDEPsiSupportRdRp}
\langle \frac{\partial U}{\partial t} + \mathcal{L}_{g}U , \psi \rangle = \langle X , \psi \rangle, \quad \forall \psi \in \SchwartzRdxR \hbox{ such that } \supp{\psi} \subset \RdxRp.
\end{equation}
\end{itemize} 
\end{Theo}

\textbf{Proof: } Let $ Y = \mathscr{F}_{S}(X) \in \mathscr{M}_{SG}(\RdxRp)$ and $ V_{0} = \mathscr{F}_{S}(U_{0}) \in \mathscr{M}_{SG}(\Rd)$. Let us then consider the solution $V$ of the transformed problem \eqref{Eq:CauchyProblemTransformedFOEMFormal} obtained from Theorem \ref{Theo:ExistenceUniquenessSolutionTransformedFOEM} and which is given by \eqref{Eq:ExpressionSolutionVTransformedFOEMugly}. Let us consider its restriction to $\RdxRp$, $\mathbf{1}_{\RdxRp}V$, which is given by \eqref{Eq:RestrictionSolutionTransformedCauchyFOEM}. We define then
\begin{equation}
\label{Eq:SolutionUCauchyProblemFOEMformal}
U = \mathscr{F}_{S}^{-1}\left( \mathbf{1}_{\RdxRp}V  \right) = \mathscr{F}_{S}^{-1}\left( e^{-tg}\left( (V_{0} - \mathcal{D}_{g}^{*}(Y)_{0})\boxtimes \mathbf{1}_{\Rp} \right) + \mathcal{D}_{g}^{*}(Y)  \right).
\end{equation}
Clearly $U \in \VRdCfdRpDual$. Since $\mathbf{1}_{\RdxRp}V$ has a càdlàg-in-time representation,  $U$ has it also (Proposition \ref{Prop:CadlagRepresentationInterchangeSpatialFourier}), and its evaluation at $t=0$ is $U_{0} = \mathscr{F}_{S}^{-1}(V_{0})$. 

Following expression \eqref{Eq:TemporalDerivvative1RdxRpV} for $\frac{\partial}{\partial t}(\mathbf{1}_{\RdxRp}V )$, and since  $\frac{\partial }{\partial t} \circ \mathscr{F}_{S}^{-1} = \mathscr{F}_{S}^{-1}  \circ \frac{\partial }{\partial t} $ over $\TemperedRdxR$, it follows that $U$ satisfies, in the sense of $\TemperedRdxR$, 
\begin{equation}
\frac{\partial U }{\partial t} + \mathcal{L}_{g}U = X + \left(U_{0} - \mathscr{F}_{S}^{-1}\left(\mathcal{D}_{g}^{*}(Y)_{0} \right) \right)\boxtimes \delta_{0_{T}}.
\end{equation}
If we restrain the space of test-functions to those in $\SchwartzRdxR$ such that their supports are included in $\RdxRp$, then we obtain condition \eqref{Eq:UsolvesPDEPsiSupportRdRp} since for such kinds of test-functions we have $\psi(\cdot , 0) = 0$ and hence $ \langle (U_{0} - \mathscr{F}_{S}^{-1}\left(\mathcal{D}_{g}^{*}(Y)_{0} \right) )\boxtimes \delta_{0_{T}} , \psi \rangle = 0$. This proves the existence of such a solution.

The uniqueness is proven in a typical manner, by supposing that there are two different solutions satisfying the conditions and then taking the difference between the solutions. It follows that such difference must be of the form $U_{H} = \mathscr{F}_{S}^{-1}( e^{-tg} S \boxtimes \mathbf{1}_{\Rp} )$ for some $S \in \MsgRd$. $U_{H}$ has a càdlàg-in-time representation which must be null. It is then immediate to conclude that $S=0$, and hence there is a unique solution satisfying the desired conditions. $\blacksquare$

The solution $U \in \VRdCfdRpDual$ satisfying \eqref{Eq:CauchyProblemFOEMFormal} in the sense of Theorem \ref{Theo:ExistenceUniquenessSolutionCuachyFOEMformal} can be then described through its càdlàg-in-time representation $(U_{t})_{t \in \Rp} \subset \VRdDual$, given by
\begin{equation}
\label{Eq:CadlagTimeRepresentationSolutionCauchyFOEM}
\langle U_{t} , \varphi \rangle = \int_{\Rd} e^{-tg(\xi)} \mathscr{F}_{S}^{-1}(\varphi)(\xi) d \mathscr{F}_{S}(U_{0})(\xi) + \int_{\Rd \times\left( 0 , t \right]} e^{-(t-s)g(\xi)}\mathscr{F}_{S}^{-1}(\varphi)(\xi)d\mathscr{F}_{S}(X)(\xi,s), 
\end{equation}
for all $\varphi \in \VRd$ and for all $t \in \Rp$. This is simply $U_{t} = \mathscr{F}_{S}^{-1}(V_{t})$, with $(V_{t})_{t \in \Rp} \subset \mathscr{M}_{SG}(\Rd)$ being the càdlàg-in-time representation of the solution to the transformed problem, given by \eqref{Eq:CadlagRepresentationSolutionTransformedCauchyProblem}.

\section{Long-term asymptotic analysis}
\label{Sec:LongTermAsymptotics}

We study right now the relation between the solution $U^{\infty} \in \VRdxRDual$ presented in Section \ref{Sec:SolutionsOverRdxR} and the solution $U$ to an associated Cauchy problem presented in Section \ref{Sec:CauchyProblem}. Precisely, we consider the case where $g$ is a continuous symbol function such that $g_{R} \geq \kappa > 0 $, and we consider $X \in \VRdxRDual \cap \VRdCfdRDual$. We denote $\mathbf{1}_{\RdxRp}X := \mathscr{F}_{S}^{-1}(  \mathbf{1}_{\RdxRp} \mathscr{F}_{S}(X)    )  \in \VRdCfdRpDual $ ($X$ is not necessarily a measure, although it behaves as a measure in time, hence we can speak about its restriction to $\RdxRp$, which we decided to define in this way). We denote $U^{\infty}$ the unique solution in $\VRdxRDual$ to equation \eqref{Eq:ThePDE}, which is given by \eqref{Eq:Uinfinity}, with càdlàg-in-time representation $(U_{t}^{\infty})_{t \in \R} \subset \VRdDual$ given by \eqref{Eq:CadlagRepresentationUinfinity}. Given an arbitrary fixed $U_{0} \in \VRdDual$, we denote $U$ the unique solution to the Cauchy problem 
\begin{equation}
\label{Eq:CauchyProblemFOEMFormalXRestricted}
\left\lbrace\begin{array}{lcr}
\dfrac{\partial U}{\partial t} + \mathcal{L}_{g}U = \mathbf{1}_{\RdxRp}X & & \hbox{over} \  \RdxRp \\
U \big|_{t=0} = U_{0} & &
\end{array}
\right.
\end{equation}
satisfying the conditions in Theorem \ref{Theo:ExistenceUniquenessSolutionCuachyFOEMformal}. The càdlàg-in-time representation of $U$, $(U_{t})_{t \in \Rp} \subset \VRdDual $ is given by \eqref{Eq:CadlagTimeRepresentationSolutionCauchyFOEM} (since the integral with respect to $\mathscr{F}_{S}(X)$ is over a subset of $\RdxRp$, there is no difference in using $\mathscr{F}_{S}(X)$ or $\mathscr{F}_{S}(\mathbf{1}_{\RdxRp}X)$). The next Theorem states that, actually, the solution $U^{\infty}$ describes how the solution $U$ behaves \textit{spatio-temporally} once the time has flown long enough.

\begin{Theo}
\label{Theo:AsymptoticConvergenceUinfty}
For every $\epsilon > 0 $ and for every $\varphi \in \VRd$ there exists $t_{\epsilon,\varphi} \in \Rp$ such that 
\begin{equation}
\label{Eq:AsymptoticConvergenceUinfty}
\left| \langle U_{t}^{\infty} - U_{t} , \phi \rangle \right| < \epsilon, \quad \forall t \geq t_{\epsilon,\varphi}, \forall \phi \in \VRd \hbox{\ translation of \ } \varphi.
\end{equation}
\end{Theo}

\textbf{Proof: } Let $\varphi \in \VRd$. Let $h \in \Rd$ and let $\phi = \tau_{h}\varphi$, where $\tau_{h}$ denotes the operator translation by $h$ ($\tau_{h}\varphi (x) = \varphi(x-h)$). We have thus $\mathscr{F}_{S}^{-1}( \phi )(\xi) = e^{ih^{T}\xi}\mathscr{F}_{S}^{-1}( \varphi )(\xi) $. Using the càdlàg-in-time representations of $U^{\infty}$ and $U$, we have
\small
\begin{equation}
\label{Eq:DifferenceUtUtinfinty}
\left. \begin{aligned}
\left| \langle U_{t} - U_{t}^{\infty} , \phi \rangle \right| &= \Bigg| \int_{\Rd} e^{-tg(\xi)}\mathscr{F}_{S}^{-1}(\phi)(\xi)d\mathscr{F}_{S}(U_{0})(\xi) + \int_{\Rd\times\left( 0 , t \right]}e^{-(t-s)g(\xi)}\mathscr{F}_{S}^{-1}(\phi)(\xi) d\mathscr{F}_{S}(X)(\xi,s) \\
&\quad\quad - \int_{\Rd\times\left( - \infty , t \right] }e^{-(t-s)g(\xi)}\mathscr{F}_{S}^{-1}(\phi)(\xi)d\mathscr{F}_{S}(X)(\xi,s)   \Bigg|   \\
&= \left| \int_{\Rd}e^{-tg(\xi)} \mathscr{F}_{S}^{-1}(\phi)(\xi)d\mathscr{F}_{S}(U_{0})(\xi) - \int_{\Rd\times\left(-\infty , 0 \right]}e^{-(t-s)g(\xi)}\mathscr{F}_{S}^{-1}(\phi)(\xi) d\mathscr{F}_{S}(X)(\xi,s) \right| \\
&\leq e^{-t\kappa}\left( \int_{\Rd}\left| \mathscr{F}_{S}^{-1}(\varphi)(\xi) \right|d\left|\mathscr{F}_{S}(U_{0})\right|(\xi) + \int_{\Rd\times\left( - \infty , 0 \right]} e^{s\kappa}\left| \mathscr{F}_{S}^{-1}(\varphi)(\xi) \right| d\left| \mathscr{F}_{S}(X) \right|(\xi,s) \right).
\end{aligned} \right.
\end{equation}
\normalsize
We set $C_{\varphi} =  \int_{\Rd}\left| \mathscr{F}_{S}^{-1}(\varphi)(\xi) \right|d\left|\mathscr{F}_{S}(U_{0})\right|(\xi) + \int_{\Rd\times\left( - \infty , 0 \right]} e^{s\kappa}\left| \mathscr{F}_{S}^{-1}(\varphi)(\xi) \right| d\left| \mathscr{F}_{S}(X) \right|(\xi,s) $, and since $\mathscr{F}_{S}(U_{0}) \in \MsgRd$, $\mathscr{F}_{S}(X) \in \MsgRdxR$, and the function $s \mapsto e^{s\kappa}$ decreases faster that any polynomial as $s \to - \infty$, it follows that $C_{\varphi} < \infty $. If $C_{\varphi} = 0$ the result is straightforward. If $C_{\varphi} > 0 $, for every $\epsilon$ we chose $t_{\epsilon,\varphi}$ big enough such that $e^{-t_{\epsilon,\varphi}\kappa} < \frac{\epsilon}{C_{\varphi}}$. It follows immediately from \eqref{Eq:DifferenceUtUtinfinty} that for every $t \geq t_{\epsilon,\varphi}$,
\begin{equation}
\left| \langle U_{t} - U_{t}^{\infty} , \phi \rangle \right| < \frac{\epsilon}{C_{\varphi}} C_{\varphi} = \epsilon. \quad \blacksquare.
\end{equation}

We remark that this notion of asymptotic convergence is quite different to the classical \textit{convergence to a steady state} widely used in evolution problems, since in those cases we usually look at for a \textit{spatial} distribution to which the solution converges in some sense. It is not clear if such a steady state exists in general, nor in which space of spatial distributions should it belong.

We finish by noticing a last simple fact: if we start with a corresponding initial condition given by the evaluation at $t=0$ of $U^{\infty}$, then the solutions $U$ and $U^{\infty}$ coincide over $\RdxRp$. This indicates that $U^{\infty}$ works as a sort of \textit{spatio-temporal fixed point} for the problem.

\begin{Prop}
\label{Prop:UinftyEqualsUifInitialCondition}
If we set the initial condition $U_{0} = U^{\infty}_{0}$ in the Cauchy problem \eqref{Eq:CauchyProblemFOEMFormalXRestricted}, then its solution $U$ and the distribution $U^{\infty}$ coincide over $\RdxRp$.
\end{Prop}

\textbf{Proof: } The càdlàg-in-time representation of $U$ is given in this case by
\begin{equation}
\langle U_{t}, \varphi \rangle = \int_{\Rd}e^{-tg(\xi)}\mathscr{F}_{S}^{-1}(\varphi)(\xi)d\mathscr{F}_{S}(U_{0}^{\infty})(\xi) + \int_{\Rd\times\left(0,t\right]}e^{-(t-s)g(\xi)}\mathscr{F}_{S}^{-1}(\varphi)(\xi)d\mathscr{F}_{S}(X)(\xi,s), 
\end{equation}
for any $\varphi \in \VRd$ and $t \in \Rp$. Since the measure $\mathscr{F}_{S}(U_{0}^{\infty})$ is described by \eqref{Eq:CadlagRepresentationUinfinity} evaluated at $t=0$, one has
\begin{equation}
\int_{\Rd}e^{-tg(\xi)}\mathscr{F}_{S}^{-1}(\varphi)(\xi)d\mathscr{F}_{S}(U_{0}^{\infty})(\xi) = \int_{\Rd\times\left( - \infty , 0 \right]} e^{-(t-s)g(\xi)}\mathscr{F}_{S}^{-1}(\varphi)(\xi)d\mathscr{F}_{S}(X)(\xi,s),
\end{equation}
from where we obtain
\begin{equation}
\langle U_{t}, \varphi \rangle =  \int_{\Rd\times\left(-\infty,t\right]}e^{-(t-s)g(\xi)}\mathscr{F}_{S}^{-1}(\varphi)(\xi)d\mathscr{F}_{S}(X)(\xi,s) = \langle U^{\infty}_{t} , \varphi \rangle.
\end{equation}
Hence the càdlàg-in-time representations of $U$ and $U^{\infty}$ coincide for $t \geq 0 $. Since càdlàg-in-time representations are unique, it must follow that $U$ and $U^{\infty}$ are equal over $\RdxRp$. $\blacksquare$

\bibliography{mibib}
\bibliographystyle{apacite}

\end{document}